\documentclass[12pt]{article}
\input{epsf}
\usepackage{amsmath}
\usepackage{amsfonts}
\usepackage{amscd}
\newtheorem{Def}{Definition}[section]
\newtheorem{Lem}[Def]{Lemma}
\newtheorem{Thm}[Def]{Theorem}

\newtheorem{Ex}[Def]{Example}

\newtheorem{Rem}[Def]{Remark}
\def\block{\hbox{${\vcenter{\vbox{\hrule height 0.4pt\hbox{\vrule 
width 0.4pt 
height 6pt \kern 5pt\vrule width 0.4pt}\hrule height 0.4pt }}}$}}
\def\qed{\hfill\block\vskip.2cm}

\def\H{{\cal H}}
\def\G{{\cal G}}
\def\A{{\cal A}}
\def\O{{\cal O}}
\def\U{{\cal U}}
\def\V{{\cal V}}

\def\N{\mathbb{N}}
\def\Z{\mathbb{Z}}

\def\R{\mathbb{R}}
\def\Co{\mathbb{C}}

\def\dim{\mbox{dim}}

\begin{document}

\title{Holonomy and parallel transport for Abelian gerbes}
\author{{\em Marco Mackaay}\footnote{presently working as a postdoc at 
the University of Nottingham, UK}\\
\'Area Departamental de Matem\'atica,\ UCEH \\ Universidade do Algarve 
\\ 8000 Faro \\ Portugal \\ e-mail: mmackaay@ualg.pt\\ \\ 
{\em Roger Picken}\\ Departamento de Matem\'{a}tica\\ Instituto 
Superior T\'{e}cnico\\ Av. Rovisco Pais\\ 1049-001 Lisboa\\ Portugal\\ 
e-mail: picken@math.ist.utl.pt}
\date{July 17, 2001}
\maketitle
\begin{abstract}
In this paper we establish a one-to-one correspondence between
$U(1)$-gerbes with connections, on the one hand, and their holonomies, for
simply-connected manifolds, or their parallel transports, in the general
case, on the other hand. This result is a higher-order analogue of the
familiar equivalence between bundles with connections and their holonomies
for connected manifolds. The holonomy of a gerbe with group $U(1)$ on a 
simply-connected manifold $M$ is a group morphism from the thin second homotopy
group to $U(1)$, satisfying a smoothness condition, where a homotopy between
maps from $[0,1]^2$ to $M$ is thin when its derivative is of rank $\leq 2$. 
For the non-simply connected case, holonomy is replaced by a parallel transport
functor between two special Lie groupoids, which we call Lie 2-groups. 
The reconstruction of the gerbe
and connection from its holonomy is carried out in detail for the 
simply-connected case. 
\end{abstract}
\section{Introduction}
In~\cite{Ba91} Barrett studied the holonomy of connections in 
principal bundles 
and proved a reconstruction theorem which showed that in a very precise 
sense all information about the connections and the bundles is contained in 
their holonomy. In this paper we obtain analogous results for Abelian gerbes 
with connections and their holonomy. 

We review the background of our work. Caetano and the second author~\cite{CP94} used a slightly different approach 
to obtain Barrett's results, which has some technical advantages. Let us 
sketch these results.  
It is well known that the holonomy of a connection in a principal $G$-bundle, 
$P$, defined over a connected smooth manifold $M$, assigns an element of 
$G$ to each smooth (based) loop in $M$. The holonomy of the composite of two 
based 
loops is exactly the product of the two holonomies. Unless the connection is 
flat, two 
homotopic loops have different holonomies in general. However, when there is a 
homotopy between the loops whose differential has rank at most 1 everywhere, 
the holonomies around the two 
loops are the same. We call these homotopies {\it thin} homotopies.  
One glance at any introductory book on algebraic topology shows that the 
homotopies used in the proof that the fundamental group obeys the group 
axioms are all thin (after smoothing at a finite number of non-differentiable 
points). Therefore the holonomy map descends to a group 
homomorphism from the {\it thin} fundamental group of $M$, $\pi^1_1(M)$, 
to $G$. We should remark here that we are using the equivalence relation 
on loops from~\cite{CP94}, while we 
are borrowing Barrett's terminology, which he used, strictly speaking, for a 
slightly different equivalence relation. In~\cite{CP94} the term 
{\it intimacy relation} was used for what we call a 
{\it thin} homotopy. In the 
present  
paper this terminological twist should not lead to any confusion, since the 
paper is intended to be self-contained. Barrett's main result, also obtained 
in the setup of~\cite{CP94}, shows that it is 
possible to reconstruct the bundle and the connection, up to equivalence, 
from the holonomy. This result is a nice strengthening of the well-known 
Ambrose-Singer theorem~\cite{AS53}.
 
Caetano and the second author~\cite{CP98} also defined the higher order thin homotopy 
groups of a manifold. The definition of the $n$th-order thin homotopy group, 
$\pi^n_n(M)$, is 
quite simple: one takes the definition of the ordinary $\pi_n(M)$, but, 
instead of dividing by all homotopies, one only divides by homotopies whose 
differentials have rank at most $n$. Just like the ordinary homotopy groups, 
all thin homotopy groups of order at least 2 are Abelian.     
Once we understand that a $G$-bundle with connection is equivalent to a 
smooth group homomorphism $\pi^1_1(M)\to G$, we can ask for a geometrical 
interpretation of a smooth group homomorphism from the second thin 
homotopy group, 
$\pi^2_2(M)$, to an Abelian group, for example the circle group, $U(1)$. As a 
matter of fact 
this was the main question left open in~\cite{CP98}.
As we show in this paper the answer is that, for a 1-connected manifold 
$M$, smooth group homomorphisms $\pi^2_2(M)\to U(1)$ correspond bijectively to 
equivalence classes of {\it $U(1)$-gerbes with connections}. 

Gerbes were first introduced by Giraud~\cite{Gi71} in an attempt to understand 
non-Abelian cohomology. Ironically only Abelian gerbes have developed into a 
nice geometric theory so far. 
Several people~\cite{Bry93,Ch98,Gaj97,Hi99} have studied 
Abelian gerbes. In this paper when we say gerbe we always mean a 
$U(1)$-gerbe. Just as a line-bundle on $M$ can be defined by a set of 
transition functions on double intersections of open sets in a covering of 
$M$, a 
gerbe can be defined by a set of ``transition line-bundles'' on double 
intersections. This 
point of view was worked out by Chatterjee~\cite{Ch98} (see also~\cite{Hi99}). The 
interest of gerbes resides in the possibility of doing differential geometry with 
them. One can define gerbe-connections and gerbe-curvatures, just as for 
bundles. Gerbe-connections and line-bundle connections have a lot in common: 
for example, the 
Kostant-Weil integrality theorem has a gerbe analogue~\cite{Bry93,Ch98}. The 
main 
difference 
is a shift in dimension: equivalence classes of line-bundles on $M$ are 
classified by $H^2(M,\mathbb{Z})$, whereas equivalence classes of gerbes on 
$M$ are classified by $H^3(M,{\mathbb{Z}})$. The curvature of an ordinary 
connection is a 2-form, the curvature of a gerbe is given by a 3-form. The 
holonomy of an ordinary connection associates a group element to each loop, 
the holonomy of a gerbe-connection associates a group element to each 2-loop, 
i.e. a smooth map from $S^2$ to the manifold. Analogously one can define 
$n$-gerbes for any $n\in\N$, which are classified by the elements of 
$H^{n+2}(M,{\mathbb{Z}})$. Note that in this convention an ordinary gerbe 
is a $1$-gerbe. The main thing to notice is the existence of a 
``geometrical ladder''~\cite{Ch98,Hi99}: a line-bundle is given by a set of 
transition functions, a gerbe is given by a set of transition line-bundles, 
a 2-gerbe is given by a set of transition gerbes, etc. Whereas two functions 
are either equal or not, two different line-bundles can also be isomorphic. 
In the definition of a gerbe it is this extra degree of freedom that 
is of most interest. On a triple intersection of open sets one requires the 
product of the transition bundles to be isomorphic to the trivial bundle by a 
given isomorphism and requires these isomorphisms to satisfy a cocycle 
condition on every four-fold intersection. The higher one gets up the 
ladder, the more intricate the notion of isomorphism becomes. One gets 
isomorphisms between isomorphisms, etc, up to the highest level where one 
requires a cocycle condition to hold. This feature is known as {\it 
categorification}, and we recommend the reader to read~\cite{BD95,BD98} on the 
general concept of categorification. Finally we note that there is 
also a notion of $n$-gerbe connection, and the parallel 
transport of an $n$-gerbe connection is defined along $n+1$-dimensional paths. 
Some details of these notions remain to be worked out completely, but 
Gajer~\cite{Gaj97,Gaj99} has worked out a considerable part already.

The truly categorical nature of gerbe-connections 
comes to light when we study the general setting of their 
parallel transport. To understand the need for this categorical language, 
one ought to think first about the right 
way of formulating parallel transport for a connection on a principal bundle 
over a manifold which is not necessarily connected. Since the concept of 
holonomy involves the choice of a base-point, one cannot expect to recover 
all information about the connection by looking only at its holonomy around 
based closed loops. One has to give up working with closed loops 
only and start working with arbitrary paths with arbitrary endpoints. 
This way one finds that the right language is that of Lie groupoids and 
groupoid morphisms, i.e. functors, instead of Lie groups and group 
homomorphisms. 
A thorough account of Lie groupoids and 
their history can be found in~\cite{Mack87}. We note that a Lie groupoid with 
only one object is precisely a Lie group. We will show how 
to translate Barrett's~\cite{Ba91} results to this more general framework. 
As pointed out by Brylinski~\cite{Bry93}, a gerbe with gerbe-connection on 
$M$ gives rise to a line-bundle with connection on $\Omega(M)$, the free 
loop space of $M$. Now the point is that, if $M$ is not 1-connected, then 
$\Omega(M)$ is not connected, so we cannot expect to recover the 
gerbe-connection just from its holonomy around based 2-loops.  
Our result about the parallel 
transport of gerbe-connections is that, if $M$ is connected at least, 
their proper setting is that of Lie groupoids with a monoidal structure 
satisfying the group laws, which we 
therefore propose to call {\it Lie 
2-groups}, and functors between them. As we show, the notion of thin 
homotopy is central in this whole story. 

We should remark that, 
when finishing this paper, we found an article by 
Gajer~\cite{Gaj99} in which he already obtains the characterization of Abelian 
$n$-gerbes with $n$-gerbe connections in terms of holonomy maps. However, 
Gajer's approach is very different from ours. Let us briefly 
explain where they differ. The main motivation for our 
approach is that we are trying to 
understand the differential geometry behind the four-dimensional state-sum 
models defined by the first author of the present paper~\cite{Ma99,Ma00}. 
For the understanding of these 
state-sums it would be helpful to find a categorical way of perceiving 
the relation between homotopy theory and differential geometry. Concretely 
the motivation was to understand what the homotopy 
$2$-type of a manifold has to do with differential geometry. We believe that 
the setup in this paper provides such a link, because in the most general 
case the parallel transport of a gerbe with gerbe-connection on $M$ 
yields a functor whose domain is the {\it thin homotopy 2-type} of $M$. 
Gajer defines everything in terms of groups and group homomorphisms and 
therefore his approach does not make this link with homotopy theory. 
For more information about our motivation see 
the beginning of Sect.~\ref{nonsimp}. 

We work out the simply-connected case ``by hand'', i.e. without relying on 
abstract cohomological arguments as Gajer does. Given the gerbe-holonomy 
we compute explicitly the gerbe and the gerbe-connection that correspond to 
it. We also remark that Gajer does not use some type of higher thin homotopy 
groups, which conceptually are easy to understand, but he uses the rather 
more complicated iterated construction, $G^{(n-1)}\left(G_{ab}(M)\right)$, 
where $G$ is a certain variant of the thin fundamental group and 
$G_{ab}$ its Abelianization. 
We hope that our approach makes the subject accessible to a 
broader group of mathematicians and physicists. 

Finally, our setup provides a link with the theory of double 
Lie groupoids  which is not 
necessarily restricted to the Abelian case~\cite{Mack92,Mack00}. We are only aware of one type of 
concrete examples of non-Abelian gerbes, which define the obstruction to 
lifting a principal $G$-bundle to a $\hat{G}$-bundle, where $\hat{G}$ is a 
non-abelian extension of $G$. To our knowledge no 
way has been found to even start a theory of connections 
on non-Abelian gerbes, and in the aforementioned concrete examples one can 
easily see that a straightforward attempt to generalize the Abelian approach 
to connections fails. Maybe our link with double Lie groupoids can shed 
some new light on non-Abelian gerbes and the possibility of defining 
connections on them. We should remark that Thompson~\cite{Th00} has 
worked out in 
detail the definition of quaternionic gerbes and some of their properties. 
Unfortunately his only examples, related to conformal $4$-manifolds, 
are actually just \v{C}ech 1-coboundaries with values in $SU(2)$ and 
therefore trivial as gerbes. For an introduction to the 
general theory of non-Abelian gerbes and its history one can 
read~\cite{Bre94}. However, the 
only concrete example in the book is the one we mentioned above.

We now give a short table of contents:
\begin{description}
\item{\ref{Int})} {\bf Introduction}.
\item{\ref{LBG})} {\bf Line-bundles and gerbes with connections}. In this 
section 
we recall some basic definitions and facts about the objects mentioned in 
the title. We claim no originality, but we hope that this section helps the 
non-specialist to understand the paper. In the same spirit we have tried to 
make the paper as self-contained as possible.   
\item{\ref{HLB})} {\bf Ordinary holonomy}. We show in detail the 
equivalence between line-bundles with connections, 
given in terms of transition 
functions and local $1$-forms, and their holonomies, as a preparation 
for the discussion of gerbe-holonomy and parallel transport.
\item{\ref{HH})} {\bf Thin higher homotopy groups}. We recall the definition 
of the thin higher homotopy groups.
\item{\ref{HG})} {\bf Gerbe-holonomy}. We show how the holonomy of a 
gerbe-connection on a gerbe on $M$ gives rise to a smooth group homomorphism 
$\pi_2^2(M)\to U(1)$. We also give a local formula for the holonomy in terms 
of the connection 0- and 1-forms, which constitute the gerbe-connection. 
\item{\ref{nonsimp})} {\bf Parallel transport in gerbes}. In this section 
we explain the key idea of the whole paper. We first show how parallel 
transport in ordinary bundles can be formulated in terms of Lie groupoids 
and smooth functors. Then we show how a gerbe with 
connection on $M$ gives rise to a Lie 2-group (a Lie 2-groupoid with only 
one object) on $M$ and how the parallel transport yields a smooth functor 
between Lie $2$-groups. 
Specialists can jump to this section right after the introduction.   
\item{\ref{BL})} {\bf Barrett's lemma for 2-loops}. This is just a 
technical intermezzo which is necessary for the sequel. On a first reading 
one can jump over this section without loss of comprehension.
\item{\ref{Rec})} {\bf The 1-connected case}. Here we reconstruct 
explicitly the \v{C}ech 2-cocycle of a gerbe and the 0- and 1-connections 
from a given holonomy map $\pi_2^2(M)\to U(1)$ for a 1-connected 
manifold $M$. In this case the reconstruction statement 
is much easier to formulate and understand than in the general case, 
because the language of groups and group homomorphisms is more 
familiar than the language of Lie 2-groups and functors between those to most 
mathematicians. We 
have therefore decided to prove this case in great detail and only indicate 
in Sect.~\ref{nonsimp} what changes have to be made to obtain the proof for 
the general case (which are only small when the change in language is well 
understood). On a first reading one might try to read this section before 
Sect.~\ref{nonsimp}.
\end{description}

\label{Int}
\section{Line bundles and gerbes with connections}
This section is just meant to recall some of the basic facts about 
line-bundles and gerbes. We claim no originality in this section. 
As a matter of fact we follow Hitchin~\cite{Hi99} and Chatterjee~\cite{Ch98} closely. 
However, we feel that a section like this is necessary, because gerbes are 
still rather unfamiliar mathematical objects to most mathematicians, and 
therefore the lack of an introductory section might scare off people who 
would like to read this paper. Of course we assume some familiarity with the 
differential geometry of principal bundles and connections.

Throughout this paper let $M$ be a smooth connected 
finite-dimensional manifold and 
$\U=\left\{U_i\ \vert\ i\in J\right\}$ an open cover of $M$ such that every 
non-empty $p$-fold intersection 
$U_{i_1\ldots i_p}=U_{i_1}\cap\cdots\cap U_{i_p}$ is 
contractible. We also assume the existence of a partition of unity 
$\left(\rho_i\colon U_i\to \R\right)$ subordinate to $\U$.
In this paper a {\it complex line-bundle} $L$ on $M$ 
can come in essentially two different but equivalent forms: as a complex 
vector bundle of rank 1 or as 
a set of transition functions $g_{ij}\colon U_{ij}\to U(1)$. 
Note that we take the values of our transition functions to be in $U(1)$ 
rather than $\Co^*$, thinking of them as transition functions of a 
principal bundle rather than a vector bundle. 
Recall that the $g_{ij}$ have to satisfy $g_{ji}= g_{ij}^{-1}$ and 
the \v{C}ech cocycle condition 
$\delta g_{ijk}=g_{ij}g_{jk}g_{ki}= 1$ on 
$U_{ijk}$. Two sets of 
transition functions $g_{ij}$ and $g'_{ij}$ define equivalent line-bundles 
if and only if there exist functions $h_{i}\colon U_i\to U(1)$ such that 
$g'_{ij}g^{-1}_{ij}= h_i^{-1}h_j$. In this case we say that 
$g_{ij}$ and $g'_{ij}$ are cohomologous and $\delta h_{ij}=h_i^{-1}h_j$ is 
the coboundary by which they differ. Thus equivalence classes of line-bundles 
correspond bijectively to cohomology classes in the first \v{C}ech 
cohomology group $\check{H}^1(M,\underline{U(1)}_M)$. Following the 
usual notation in the literature on gerbes, we write $\underline{G}_M$ 
for the sheaf of smooth $G$-valued functions on (open sets of) $M$, 
where $G$ is any Lie group. Following the same convention, the sheaf of 
{\it constant} $G$-valued functions is denoted by $G_M$.     
The bijection above defines a group homomorphism, where 
the group operation for cochains is defined by pointwise multiplication and 
for line-bundles by their tensor product. Similarly we use two different but 
equivalent ways to write down a 
{\it connection} in $L$: as a 
covariant derivative $\nabla\colon\Gamma(M,L)\to\Omega^1(M,\Co)$ or 
as a set of of local 1-forms $A_i\in\Omega^1(U_i)$. For each $i$ the 1-form 
$A_i$ is of course the pull-back of the connection 1-form associated to 
$\nabla$ via a local section $\sigma_i\colon U_i\to p^{-1}(U_i)\subset L$. 
We will not say more about 
covariant derivatives, but we recall that the $A_i$ have to satisfy the rule 
$$i(A_j-A_i)=d\log g_{ij}.$$  
Two line-bundles, $g_{ij}$ and $g'_{ij}$, with connections, $A_i$ and 
$A'_i$ respectively, are equivalent exactly when $g_{ij}$ and $g'_{ij}$ 
define equivalent bundles 
as above and $iA'_i= iA_i+d\log h_i$. 
Given a line-bundle, $g_{ij}$, with connection, $A_i$, we can define 
its {\it curvature 2-form} 
$F\in\Omega^2(M)$ by $F\vert_{U_i}=dA_i$. Using the partition of unity 
$\left(\rho_i\right)$ we can easily define a connection for a given 
line-bundle, $g_{ij}$, by $A_i=-i\sum_{\alpha}\rho_{\alpha}d\log g_{\alpha i}$. 
An immediate consequence of the definitions is that the cohomology class of 
$F$ is independent of the chosen connection. The invariant $[F]/2\pi$ 
is called the Chern class of the line-bundle. A well known fact about 
line-bundles is the Kostant-Weil integrality theorem (see~\cite{Bry93} and 
references therein), which says that 
any closed 2-form $F$ on $M$ is the curvature 2-form of a connection in a 
line-bundle if and 
only if $[F]/2\pi$ is the image of an 
integer singular cohomology 
class in $H^2(M,\Z)$. Basically this theorem is the statement that the 
cohomology groups $\check{H}^1(M,\underline{U(1)}_M)$ and $H^2(M,\Z)$ are 
isomorphic, which 
follows from the exact sequence of sheaves 
$$0\longrightarrow\Z_M\stackrel{\times 2\pi i}
{\longrightarrow}\underline{\Co}_M\stackrel{\exp}{\longrightarrow}
\underline{\Co}^*_M\longrightarrow 1$$
and from the isomorphism 
$$\check{H}^2(M,\Z_M)\cong H^2(M,\Z).$$
  
Finally we should 
remark 
that Deligne (see references in~\cite{Bry93}) 
found a nice way to encode line-bundles with 
connections in a cohomological terminology. Equivalence classes of 
line-bundles with connections correspond bijectively to cohomology classes 
in the first (smooth) Deligne hypercohomology group 
$H^1(M,\underline{\Co}^*\stackrel{d\log}{\rightarrow}\underline{\A}^1_{M,\Co})$. 
This means little more than that a line-bundle with connection can be 
defined by a pair 
of local data $\left(g_{ij},A_i\right)$ with respect to the covering 
$\U$ such that 
$$D\left(g_{ij},A_i\right)=
\left(\delta g_{ijk},i(\delta A_{ij}-d\log g_{ij})\right)=(1,0)$$ 
and that two such pairs $\left(g_{ij},A_i\right),
\left(g'_{ij},A'_i\right)$ represent equivalent line-bundles with 
equivalent connections precisely if there is a set of functions 
$h_i$ such that 
$$\left(g'_{ij}g_{ij}^{-1},i(A'_i-A_i)\right)=
\left(\delta h_{ij},d\log h_i\right)=D(h_i).$$ 
The difference between this remark and the precise definition of Deligne 
cohomology lies in the fact that the local data are not just a family of 
sets but a {\it sheaf}. This distinction is important for a rigorous 
treatment because one wants all definitions to be independent of the choice 
of open cover in the end. However, this is not the right place to define the 
whole machinery of sheaves and sheaf cohomology. For a rigorous introduction 
to Deligne hypercohomology see Brylinski's book~\cite{Bry93}.

Just as in the case of line-bundles we can define a gerbe by two different but 
equivalent kinds of data. The first alternative is to say that a gerbe is 
given by a set of {\it 
transition line-bundles} $\Lambda_{ij}$ on $U_{ij}$ together with a 
nowhere zero section (also called trivialization) 
$\theta_{ijk}\in\Gamma(U_{ijk},\Lambda_{ijk})$ of the tensor product 
line-bundle $\Lambda_{ijk}=\Lambda_{ij}\Lambda_{jk}\Lambda_{ki}$. The 
transition line-bundles have to satisfy $\Lambda_{ji}\cong\Lambda_{ij}^{-1}$, 
where the latter is the inverse of $\Lambda_{ij}$ with respect to the 
tensor product, and the cocycle condition which says that 
$$\Lambda_{ijk}=\Lambda_{ij}\Lambda_{jk}\Lambda_{ki}\cong 1,$$ 
where the last line-bundle is the trivial line-bundle over $U_{ijk}$. The 
trivializations have to satisfy $\theta_{s(i)s(j)s(k)}=
\theta_{ijk}^{\epsilon(s)}$, for any permutation $s\in S_3$, where 
$\epsilon(s)$ is the sign of $s$, and the cocycle 
condition $$\delta\theta_{ijkl}=\theta_{jkl}\theta_{ikl}^{-1}\theta_{ijl}
\theta_{ijk}^{-1}=1.$$
To understand the last equation one has to note that on $U_{ijkl}$ the 
tensor product of all the line-bundles involved is identical to the trivial 
line-bundle, of course up to the fixed isomorphisms $\Lambda_{ji}\cong
\Lambda_{ij}^{-1}$ and $\Lambda_{ij}\Lambda_{ij}^{-1}\cong 1$, and the 
canonical 
isomorphism which reorders the factors in the tensor product. This 
holds because every line-bundle appears twice in the product with 
opposite signs. The last equation means that the product of the sections has 
to be equal, forgetting about the uninteresting isomorphisms above, to the 
canonical section in the trivial line-bundle over $U_{ijkl}$ (i.e. 
the one which associates to each $x\in U_{ijkl}$ the point $(x,1)\in 
U_{ijkl}\times U(1)$).

Two gerbes are said to be equivalent if for each $i\in J$ there exists a 
line-bundle $L_i$ on $U_i$, such that for each $i,j\in J$ we have
bundle isomorphisms
$$m_{ij}\colon L_j\stackrel{\cong}{\rightarrow}\Lambda'_{ij}
\Lambda_{ij}^{-1}\otimes L_i$$
such that 
$$m_{ij}\circ m_{jk}\circ m_{ki}=\theta'_{ijk}\theta_{ijk}^{-1}\otimes 
\mbox{id}$$
on $U_{ijk}$. 
The data defining an equivalence are called an {\it object} by 
Chatterjee~\cite{Ch98}. 

In order to relate the definition above to \v{C}ech 
cohomology we only have to remember that all our $p$-fold intersections are 
contractible so all line-bundles above are necessarily equivalent to the 
trivial line-bundle. This means that we can choose a nowhere zero section 
$\sigma_{ij}$ in $\Lambda_{ij}$ for each $i,j\in J$. With respect to these 
sections we now get 
$\theta_{ijk}=g_{ijk}\sigma_{ijk}$, where we take 
$\sigma_{ijk}=\sigma_{ij}\sigma_{jk}\sigma_{ki}$ in the tensor product 
line-bundle. By the definition above we see that $g$ defines a \v{C}ech 
2-cocycle, i.e., for all $i,j,k,l\in J$ we have 
$$\delta g_{ijkl}=g_{jkl}g^{-1}_{ikl}g_{ijl}g^{-1}_{ijk}=1$$
on $U_{ijkl}$. Given two equivalent gerbes $(\Lambda'_{ij},\theta'_{ijk}), 
(\Lambda_{ij},\theta'_{ijk})$ and an object $(L_i,m_{ij})$ for 
$\Lambda_{ij}'\Lambda_{ij}^{-1}$, we can 
also choose 
nowhere zero sections $\sigma_i$ in $L_i$ for each $i\in J$. This gives us 
two nowhere zero sections in $\Lambda'_{ij}\Lambda_{ij}^{-1}
\otimes L_j$: $\sigma'_{ij}\sigma_{ij}^{-1}\sigma_j$ 
and $m_{ij}\circ\sigma_i$. For each $i,j\in J$ we can take the quotient of 
these sections which defines a function $h_{ij}\colon 
U_{ij}\to U(1)$ such that for each $i,j,k\in J$ we have 
$$g'_{ijk}g^{-1}_{ijk}=\delta h_{ijk}=h_{jk}h_{ik}^{-1}h_{ij}$$
on $U_{ijk}$. Thus every equivalence class of gerbes induces a \v{C}ech 
cohomology class in $\check{H}^2(M,\underline{U(1)}_M)$. Hitchin~\cite{Hi99} shows that the 
converse is true as well, which leads to the conclusion that equivalence 
classes of gerbes correspond bijectively to the elements in 
$\check{H}^2(M,\underline{U(1)}_M)$.

Gerbe-connections are defined by two sets of data: the 0-connections and the 
1-connections (Chatterjee's terminology). Let $\G$ be a gerbe given by 
a set of line-bundles $\Lambda_{ij}$ and trivializations $\theta_{ijk}$. A 
{\it 0-connection} in $\G$ consists of a covariant derivative $\nabla_{ij}$ 
in $\Lambda_{ij}$ for each $i,j\in J$ such that for each $i,j,k\in J$
$$\nabla_{ijk}\theta_{ijk}=0$$
on $U_{ijk}$.
Here $\nabla_{ijk}$ is the covariant derivative in the tensor product of the 
corresponding line-bundles induced by the covariant derivatives in these 
line-bundles. Given two line-bundles , $L$ and $L'$, with connections, 
$\nabla$ and $\nabla'$ respectively, Brylinski~\cite{Bry93} denotes the 
induced connection in 
$L\otimes L'$ by $\nabla+\nabla'$, which is defined by 
$$\left(\nabla+\nabla'\right)\sigma\otimes\sigma'=
\nabla\sigma\otimes\sigma'+\sigma\otimes\nabla'\sigma',$$
for any sections $\sigma\in L$ and $\sigma'\in L'$. We also require that 
$\nabla_{ij}+\nabla_{ji}=0$ for all $i,j\in J$. 
The alternative definition of a 
0-connection is obtained in a straightforward manner by taking the pull-back 
of the connection 1-form associated to $\nabla_{ij}$ via $\sigma_{ij}$ for 
each $i,j\in J$. Let $\G$ correspond to the 2-cocycle 
$g_{ijk}$, and choose a logarithm of $g_{ijk}$. 
A 0-connection can then be defined by a set of 1-forms 
$A_{ij}\in\Omega^1(U_{ij})$ 
such that 
$$i(A_{ij}+A_{jk}+A_{ki})=-d\log g_{ijk}.$$ 
Of course we assume that $A_{ji}=-A_{ij}$. 
A {\it 1-connection} in $\G$ is defined by a set of local 2-forms 
$F_i\in\Omega^2(U_i)$ such that 
$$F_j-F_i=\sigma^*_{ij}K(\nabla_{ij}),$$ 
where $K(\nabla_{ij})$ denotes the curvature of $\nabla_{ij}$. Alternatively 
we get 
$$F_j-F_i=dA_{ij}.$$
In the sequel we sometimes denote $dA_{ij}$ by $F_{ij}$. A 0-connection and a 
1-connection 
on $\G$ together form what we call a {\it gerbe-connection}. Our typical 
notation for a gerbe-connection is $\A$. In a natural way a 
gerbe-connection, $\A$, leads to the notion of a gerbe-curvature 
3-form, $G\in\Omega^3(M)$, 
which is defined by $G\vert_{U_i}=dF_i$. Again using the partition of unity 
we see that, for a given gerbe $g_{ijk}$, it is easy to 
define a 0-connection by 
$A_{ij}=-i\sum_{\alpha}\rho_{\alpha}d\log g_{\alpha ij}$ and 
a 1-connection by $F_i=\sum_{\beta}\rho_{\beta}dA_{\beta i}$. 
Just as for line-bundles one can show that the cohomology class of the 
gerbe-curvature, $[G]$, does not depend on the chosen 0- and 1-connection 
(the proof is a bit harder though), and that any closed 3-form $G$ on $M$ is a 
gerbe-curvature 3-form if and only if 
$[G]/2\pi$ is the image of a class in $H^3(M,\Z)$. The cohomology 
class $[G]/2\pi$ is called the Dixmier-Douady class of the gerbe.

To complete the picture we have to 
define when two gerbes, $\G$ and $\G'$, with gerbe-connections, $\A$ and 
$\A'$ respectively, are 
equivalent. First of all, such an equivalence requires $\G$ and $\G'$ to be 
equivalent as gerbes. Let $\left(L_i,m_{ij}\right)$ define an object for this 
equivalence and let $h_{ij}$ be the \v{C}ech 
cochain corresponding to this object. $\A$ and $\A'$ are now equivalent if 
for every $i\in J$ there 
exists a connection $\nabla_i$ in $L_i$ such that $m_{ij}$ maps 
\begin{equation}
\nabla_j\mapsto\nabla'_{ij}-\nabla_{ij}+\nabla_i
\label{oc0}
\end{equation}
and such that 
\begin{equation}
F'_i=F_i+\sigma^*_i K(\nabla_i).
\label{oc1}
\end{equation} 
Equivalently this means that for every 
$i\in J$ there 
exists a 1-form $A_i\in\Omega^1(U_i)$ such that 
\begin{equation}
iA'_{ij}=i(A_{ij}+A_j-A_i)-d\log h_{ij}
\label{loc0}
\end{equation}
on $U_{ij}$ and 
\begin{equation}
F'_i=F_i+dA_i
\label{loc1}
\end{equation}
on $U_i$. Now there is a subtlety in Chatterjee's~\cite{Ch98} 
terminology. Suppose a gerbe $\G$ can be trivialized by an object 
$(L_i,m_{ij})$ (or, equivalently, by $h_{ij}$). He calls local data $\nabla_i$ 
(or, equivalently, $A_i$) which satisfy (\ref{oc0}) (resp. (\ref{loc0})) 
an {\it object $0$-connection}. 
In general such an object $0$-connection need not satisfy 
(\ref{oc1}) (resp. (\ref{loc1})). 
That would only be possible if the object-connection were 
trivializable. As Chatterjee proves, on a $2$-dimensional surface 
any gerbe with any gerbe-connection 
admits an object with an object $0$-connection, but in general the latter 
does not satisfy (\ref{oc1}). This is of course analogous to the fact that 
any line-bundle on the circle admits a trivialization, but an 
arbitrary connection in such a line-bundle need not be trivializable. 
In this article we call an object $0$-connection simply an object connection. 
It is important to keep these remarks in mind for Sect.~\ref{HG}. 

Finally we should remark that there is a bijective 
correspondence between equivalence classes of gerbes with 
gerbe-connections and cohomology classes in $H^2(M,\underline{\Co}^*\stackrel{d\log}{\rightarrow}\underline{\A}^1_{M,\Co}\stackrel{d}{\rightarrow}\underline{\A}^2_{M,\Co})$, the 
next order Deligne hypercohomology group. Again this amounts to little more 
than saying that a gerbe with gerbe-connection is defined by a triple of 
local data $\left(g_{ijk},A_{ij},F_i\right)$ satisfying the conditions 
we have explained above.

\begin{Ex}{\rm Let $G$ be a Lie group and $1\rightarrow U(1)\stackrel{i}
{\rightarrow}\hat{G}\stackrel{\pi}{\rightarrow}G\rightarrow 1$ a central 
extension. It is well known that any central extension is locally 
split. This means that $\hat{G}\stackrel{\pi}{\rightarrow}G$ is a locally 
trivial principal $U(1)$-bundle. Given a principal $G$-bundle over $M$, 
denoted by $P$, in the form of its transition functions 
$g_{ij}\colon U_{ij}\to G$, we can 
locally lift these transition functions to obtain 
$\hat{g}_{ij}\colon U_{ij}\to\hat{G}$ (by assuming that the image of $g_{ij}$ 
is sufficiently small so that the central extension can be trivialized 
over it). In general the $\hat{g}_{ij}$ do not define a cocycle, but 
$\pi(\delta \hat{g}_{ijk})=\delta g_{ijk}=1$ of course, so we have  
$\delta\hat{g}_{ijk}\colon U_{ijk}\to\ker{\pi}\cong U(1)$. As a matter of 
fact $\delta\hat{g}_{ijk}$ defines a gerbe, because 
$\delta^2=1$ always. Thus the obstruction to lifting 
$P$ to a principal $\hat{G}$-bundle defines a gerbe.} 
\end{Ex}
 
\begin{Ex}{\rm Let $S^3\subset\R^4$ be the three-dimensional sphere and take 
$N=S^3-\{(0,0,0,1)\}$ and $S=S^3-\{(0,0,0,-1)\}$. The intersection 
$N\cap S$ is homotopy equivalent to $S^2$. For example, 
$(x,y,w,z)\mapsto (x,y,w)/\sqrt{x^2+y^2+w^2}$ 
defines a homotopy equivalence. Therefore 
equivalence classes of line-bundles on $N\cap S$ are determined by cohomology 
classes in $H^2(S^2,\Z)=\Z$. Since there are no 3-fold intersections, any 
line-bundle on $S^2$ defines a gerbe on $S^3$. Let $\Lambda_{NS}$ be such a 
line-bundle. Any connection $A_{NS}$  in $\Lambda_{NS}$ defines a 
0-connection in the 
gerbe. A 1-connection is also easy to obtain, because the curvature 2-form 
of $A_{NS}$, which we denote by $F_{NS}$, can always be extended to $S$, 
since $S$ is contractible, and therefore we can define $F_S$ to be this 
extension of $F_{NS}$ and define $F_N$ to be zero. Chatterjee~\cite{Ch98} and 
Hitchin~\cite{Hi99} show how to obtain a gerbe for any codimension-3 
submanifold 
and how to define a gerbe-connection for such a gerbe. There is a little 
subtlety that we should explain: in this example $U_N\cap U_S$ is not 
contractible. 
However our definition of a gerbe in terms of line-bundles does not use 
that condition at all. Only when one wants to pass to the corresponding 
\v{C}ech-cocycle and the local forms which define the gerbe-connection 
one has to assume that all intersections are contractible. 
If one feels happier with contractible intersections one can subdivide $N$ and 
$S$, but of course this makes the definition of the gerbe a little bit more 
complicated. For a detailed treatment see~\cite{Ch98}.}          
\end{Ex}

\label{LBG}
\section{Ordinary holonomy}
It is well-known that a principal $G$-bundle with connection over $M$ 
allows one to define the notion of holonomy around any smooth closed
curve on $M$ (Kobayashi and Nomizu \cite{KN63}). In particular, given a fixed
basepoint $\ast$ in $M$, and a point in the fibre over $\ast$, this data induces
an assignment of an element of $G$ to each loop in $M$, based at $\ast$. Such an
assignment is called a holonomy map, or simply a holonomy.  In Barrett \cite{Ba91},
and in a slightly different fashion in Caetano-Picken \cite{CP94}, it was shown that
suitably defined holonomy maps are in one-to-one correspondence with $G$-bundles
plus connection plus the choice of a point in the fibre over $\ast$, up to
isomorphism. This result should be seen as a geometric version of the well-known
equivalence between flat $G$-bundles modulo gauge transformations over $M$ and
${\rm Hom} (\pi_1(M), G)/G$, the group homomorphisms from $\pi_1(M)$ to $G$,
modulo conjugation by elements of $G$.

The reconstruction of the bundle and connection from a holonomy map in Barrett~\cite{Ba91} and Caetano-Picken~\cite{CP94} was carried out in the total space of the
bundle, using Ehresmann connections. The main result of this section is to prove
the same equivalence using instead the local data defining a bundle and
connection from section \ref{LBG}. Since the aim is to prepare the ground for the
gerbe discussion in sections \ref{HG} and \ref{Rec}, we will only show the result
for the case $G=U(1)$.

We start by recalling briefly the definition of holonomy map from \cite{CP94}. 
Let $\Omega^\infty(M)$ be the space of smooth
loops $\ell:  [0,1]  \rightarrow M$ based at $\ast$ such that $\ell(t)=\ast$ for $0\leq
t < \epsilon$ and $1-\epsilon<t\leq 1$ for some $0<\epsilon<1/2$. We say that
the loop {\it sits}, or has a {\it sitting instant}, at $t=0$ and $t=1$. 
In~\cite{CP94} it is shown how to reparametrize any path to sit at its 
endpoints using 
a smoothly increasing function $\beta\colon[0,1]\to[0,1]$ with 
$\beta(t)=0$, for $t\in[0,\frac{1}{3}]$, and $\beta(t)=1$, for $t\in
[\frac{2}{3},1]$. As usual we may
define the product and inverse of loops, and the sitting property means that the
product closes in  $\Omega^\infty(M)$. We will denote the product of loops by $\star$. 
An equivalence relation between loops
appropriate for parallel transport purposes is the following: 
\begin{Def}
\label{rank1hom} Two loops $\ell$ and $\ell'$ belonging to $\Omega^\infty(M)$ are said to be
{\rm rank-$1$ homotopic}, or {\rm thin homotopic}, written $\ell\stackrel{1}{\sim} \ell'$, if there exists a map $H:
[0,1]\times [0,1]  \rightarrow M$ such that:
\begin{itemize}
\item[1.] $H$ is smooth throughout $ [0,1]\times [0,1]$
\item[2.] ${\rm rank} (DH_{(s,t)})\leq 1$  $ \forall (s,t) \in [0,1]\times [0,1]$
\item[3.] there exists $0<\epsilon<1/2$ such that
\begin{eqnarray*}
0\leq s \leq \epsilon & \Rightarrow & H(s,t)=\ell(t)\\
1-\epsilon \leq s \leq 1 & \Rightarrow & H(s,t)=\ell'(t)\\
0\leq t \leq \epsilon & \Rightarrow & H(s,t)=\ast\\
1-\epsilon \leq t \leq 1 & \Rightarrow & H(s,t)=\ast
\end{eqnarray*}
\end{itemize}
\end{Def}
Now the space of equivalence classes $\pi_1^1(M, \ast) = \Omega^\infty(M)/
\stackrel{1}{\sim}$ acquires the structure of a group in exactly the same way as
$\pi_1(M)$ does, but using rank-$1$ homotopy instead of ordinary homotopy, since
all homotopies used in the construction of $\pi_1(M, \ast)$ are in fact rank-$1$. Again the function $\beta$ is used to give the usual 
homotopies sitting endpoints. In the rest of this paper all paths and 
homotopies are understood to have sitting endpoints, which can always be 
achieved by reparametrization with $\beta$, as shown in~\cite{CP94}. 
We remark that in \cite{CP94} the terms {\it intimate} and the {\it 
group of loops}, ${\cal GL}^\infty (M)$, were used instead of 
{\it rank-$1$ homotopic} and $\pi_1^1(M, \ast)$ 
respectively. When the basepoint $\ast$ is understood we will frequently write 
$\pi_1^1(M)$ instead of $\pi_1^1(M, \ast)$.

In an analogous fashion we may introduce the smooth path groupoid $P_1^1(M)$, 
consisting of
smooth paths $p:I\rightarrow M$ which are constant in a neighbourhood of $t=0$ and $t=1$, 
identified up to rank-$1$ homotopy, defined as in Def.~\ref{rank1hom} with the obvious 
modification $0\leq t \leq \epsilon  \Rightarrow  H(s,t)=p(0)$, 
$1-\epsilon \leq t \leq 1  \Rightarrow  H(s,t)=p(1)$. Multiplication in 
$P_1^1(M)$, when possible, will
also be denoted by $\star$. The set of paths without dividing by the thin 
homotopy relation we denote by $P^{\infty}(M)$. 

We define a {\it smooth 
family of loops} to be a map $\psi\colon U\subseteq\R^k\to\Omega^{\infty}(M)$ 
defined on an open set $U\subseteq\R^k$ such that the function 
$\psi(x,t)=\psi(x)(t)$ is smooth on $U\times I$.

\begin{Def}
\label{smoothfam}
A {\rm holonomy} is a group morphism ${\cal H}: \pi_1^1(M)\rightarrow U(1)$,
which is smooth in the following sense: for every smooth family of loops
$\psi:U\subseteq {\bf R}^k \rightarrow \Omega^\infty (M)$, the 
composite 
$$
U\stackrel{\psi}{\rightarrow} \Omega^\infty(M) \stackrel{\rm proj}{\rightarrow} \pi_1^1(M)
\stackrel{\cal H}{\rightarrow} U(1),
$$
where $\rm proj$ is the natural projection, is smooth throughout $U$.
\end{Def}

For later it will be useful to have the following result.
\begin{Lem} (Recentering a holonomy) Let ${\cal H}: \pi_1^1(M, \ast) \rightarrow U(1)$ be a 
holonomy, $m\in M$ and $p\in P^\infty(M)$ a path from $\ast$ to $m$. Then 
$\tilde{{\cal H}}: \pi_1^1(M, m) \rightarrow U(1)$ defined by 
$\tilde{{\cal H}}(\ell)={\cal H}(p\star \ell \star p^{-1})$ is a holonomy.
\end{Lem}
{\bf Proof}: Since the smoothness property is clear, we only have to show that 
$\tilde{{\cal H}}$ is a group morphism, which follows from the fact that products
$p^{-1}\star p$ may be cancelled up to rank-$1$ homotopy. \qed

Suppose now that a $U(1)$ bundle $L$ with connection $\cal A$ is given in terms
of local data, $g_{ij}$ and $A_i$, as in section \ref{LBG}.  
Define their holonomy map ${\cal H}$ from $\Omega^\infty(M)$ to $U(1)$ as
follows:
$$
{\cal H}(\ell) = \exp i \int_I A^{\ell}
$$
where $A^{\ell}$ is a $1$-form on $I$ defined on each open set of the pullback 
cover via $\ell$ of the interval, ${\cal V}=\left\{V_i,i\in J\right\}$, by
$$
i A^{\ell} = i\ell^\ast (A_i) - d\log k_i
$$
where $k_i: V_i\rightarrow U(1)$ is an arbitrary trivialization of
$\ell^\ast L$, i.e. $\ell^*(g_{ij})=k_jk_i^{-1}$ on $V_{ij}$. 
The $1$-form $A^{\ell}$ is globally defined on $I$ since
\begin{eqnarray*}
i\ell^\ast(A_j-A_i)&=&\ell^\ast(g_{ij}^{-1}dg_{ij})\\
&=&(k_jk_i^{-1})^{-1}d(k_jk_i^{-1})\\
&=&d\log k_j - d\log k_i.
\end{eqnarray*}
Furthermore ${\cal
H}(\ell) $ doesn't depend on the choice of trivialization $k_i$,
since, if ${A'}^{\ell}=i\ell^\ast (A_i) - {k'_i}^{-1}dk'_i$, then $k_i/k'_i=k_j/k'_j$ on
$V_{ij}$, so that $f$, defined locally by $f=k_i/k'_i$, is a function on
$I$. Therefore
$$
\exp i\int_I (A^{\ell}- {A'}^{\ell}) = \exp \int_I f^{-1} df = f(1) f(0)^{-1} =1.$$

\begin{Lem}
${\cal H}$ descends to $\pi_1^1(M)$ and defines a holonomy, which is
independent of the choice of data $L, \cal A$ up to equivalence.
\end{Lem}
{\bf Proof}: Suppose that $\ell \stackrel{1}{\sim} \ell'$, and that $H:I^2\rightarrow M$
is a rank-$1$ homotopy between $\ell$ and $\ell'$, as in Def.~\ref{rank1hom}. Let
$k_i$ be a trivialization of $H^\ast(L)$ over $I^2$, and $A^H$ be the $1$-form
on $I^2$ defined locally on $H^{-1}(U_i)$ by $iA^H=iH^\ast(A_i) - d\log k_i$.
Now, 
\begin{eqnarray*}
\exp i \int_I (A^{\ell}-A^{\ell'})&=&\exp i\int_{I^2} dA^H\\
&=&\exp i \int_{I^2} H^\ast (F)\\
&=&1,
\end{eqnarray*}
using Stokes' theorem 
in the first equality and the fact that
$H$ is of rank $\leq 1$, whereas $F$ is a $2$-form, in the final equality.
Thus $\H$ descends to $\pi_1^1(M)$. 

Also ${\cal H}$ is a group homomorphism from $\pi_1^1(M)$ to $U(1)$ since 
$$
\int_IA^{\ell\star \ell'}= \int_I A^{\ell} + \int_I A^{\ell'}.
$$
Suppose $\psi:U\subset {\bf R}^k\rightarrow \Omega(M)$ is a smooth family of
loops in the sense of Def.~\ref{smoothfam}. Without
loss of generality we suppose that $U$ is contractible, and let $k_i$ be a
trivialization of the pull-back bundle ${\psi}^\ast(L)$. Once again we 
define a $1$-form on 
$U\times I$
by $iA^{\psi}= i{\psi}^\ast(A_i) - d\log k_i$ on each open set of the pullback cover
under $\psi$. Now
$$
{\cal H}(\psi(s_1, \ldots, s_k)) = 
(\exp i \int_{I} A^{\psi}) (s_1,\ldots, s_k)
$$
is a smooth function of $s_1,\ldots, s_k$, since all functions are smooth and
integration is a smooth operation. Thus $\H$ defines a holonomy. 

Finally ${\cal H}$ does not
depend on the choice of ${L,{\cal A}}$ up to equivalence, which follows 
immediately from Stokes' theorem and the local formula for holonomy which 
we give below.

\qed

\begin{figure}
\centerline{
\epsfbox{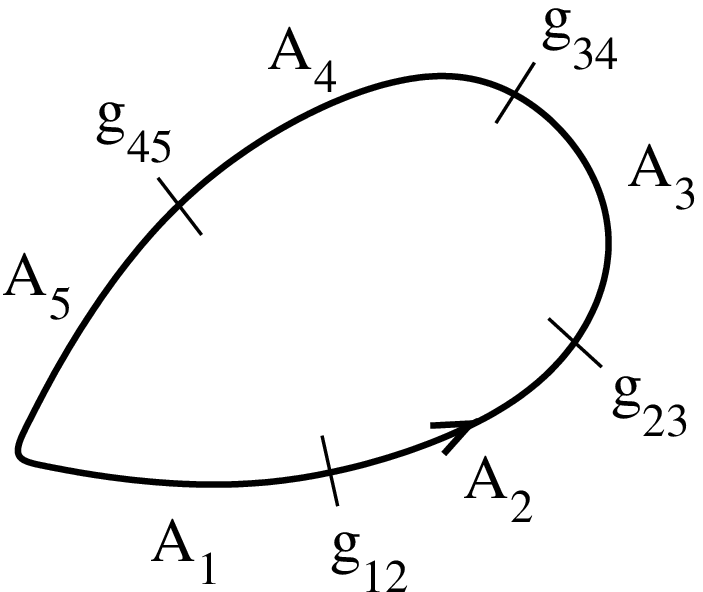}
}
\caption{local formula for holonomy}
\label{lochol1}
\end{figure}

For later purposes it is convenient to have a local expression for the holonomy
${\cal H}$, defined directly in terms of $g_{ij}$ and $A_i$. Let
$\ell:I\rightarrow M$ be a smooth loop, based at $\ast$. Fix an 
element of the open cover, $U_0$, such that $\ast \in U_0$. Consider again the 
pull-back cover
${\cal V}$ of the interval $I$. Since $I$ is a compact metric space 
it has a
Lebesgue number $\lambda>0$ such that $\forall t\in I$ we have $
]t-\lambda,t+\lambda [ \subseteq V_i$ for some $i$. Thus, given a 
decomposition of
the unit interval  $0=x_0<x_1< \cdots < x_N=1$ such that $x_\alpha - x_{\alpha-1} <\lambda, \,\forall 
\alpha=1, \ldots , N$, each subinterval $I_\alpha = [ x_{\alpha - 1},
x_{\alpha}]$  is contained in $V_{\alpha}$ for some $\alpha$, and
furthermore, by choosing a smaller $\lambda$ if necessary, or adjusting the
decomposition, we can ensure
$V_1=V_N=V_0$. Since $$
\exp \int_{I_\alpha}\left(i\ell^\ast(A_{\alpha}) - d\log k_{\alpha}\right)=
k_{\alpha}(x_{\alpha-1}) \left(\exp  \int_{I_\alpha} i\ell^\ast
(A_{\alpha})\right) k_{\alpha}^{-1}(
x_{\alpha})
$$
and $k_{\alpha}^{-1} k_{{\alpha +1}}=
\ell^\ast g_{\alpha,{\alpha + 1}}$ on $V_{\alpha,\alpha +1}$, we arrive at the 
following local formula
\begin{equation}
{\cal H}(\ell)=
\prod_{\alpha=1}^{N}\exp i\int_{I_{\alpha}}\ell^\ast(A_{\alpha}) 
\cdot g_{\alpha,\alpha+1}(\ell(x_{\alpha})),
\label{eq:locforhol}
\end{equation}
where we set $U_{N+1}=U_1$, and thus $g_{N,N+1}(\ell(x_N))=1$. 
Fig.~\ref{lochol1} sketches this local formula. 
We may also extend this
formula by an identical procedure to smooth paths $p:I\rightarrow M$. 
However, for a path that is not closed we can not in general choose 
the initial and final open sets to be equal and the formula depends on the 
specific choice (but not on the intermediate open sets in the covering of the 
path). Thus we define 
\begin{equation*}
{\cal H}^{1,N}(p)=
\left(\prod_{\alpha=1}^{N-1}\exp i\int_{I_{\alpha}}p^\ast(A_{\alpha})\cdot  
g_{\alpha,\alpha+1}(p(x_{\alpha}))\right)\cdot \exp i \int_{I_N}p^{\ast}(A_N),
\label{eq:locforpar}
\end{equation*}
and we 
have the multiplicative formula:
\begin{equation}
{\cal H}^{ij}(p\star q) = {\cal H}^{ik}(p)\cdot
g_{kl}(p(1))\cdot 
{\cal H}^{lj}(q).
\label{eq:mult}
\end{equation}

Now we turn to the reconstruction of a bundle with connection from a given holonomy
$\cal H$. Assume that the cover $\U$ is such that for each $i$ we have a 
diffeomorphism
$\phi_i:U_i\rightarrow B(0,1)$, where $B(0,1)$ is the open unit ball in ${\bf
R}^n$.  Since $M$ is path-connected we may choose a smooth path $p_i\in P^\infty
(M)$ from $\ast$ to $x_i=\phi^{-1}(0)$, the centre of $U_i$. For $U_0$ we set $p_0 $ to be the  constant path at $\ast$. Given $x\in U_i$ there is a natural path $\gamma_{i,x}\in P^\infty(M)$ 
from $x_i$
to $x$, namely the pullback under $\phi_i$ of the radial path from the 
origin to
$\phi_i(x)$ in ${\bf R}^n$, reparametrized to be constant in a neighbourhood of
$t=0,1$. For future use it is practical to define 
$$p_{i,x}=p_i\star\gamma_{i,x}.$$ 

\begin{figure}
\centerline{
\epsfbox{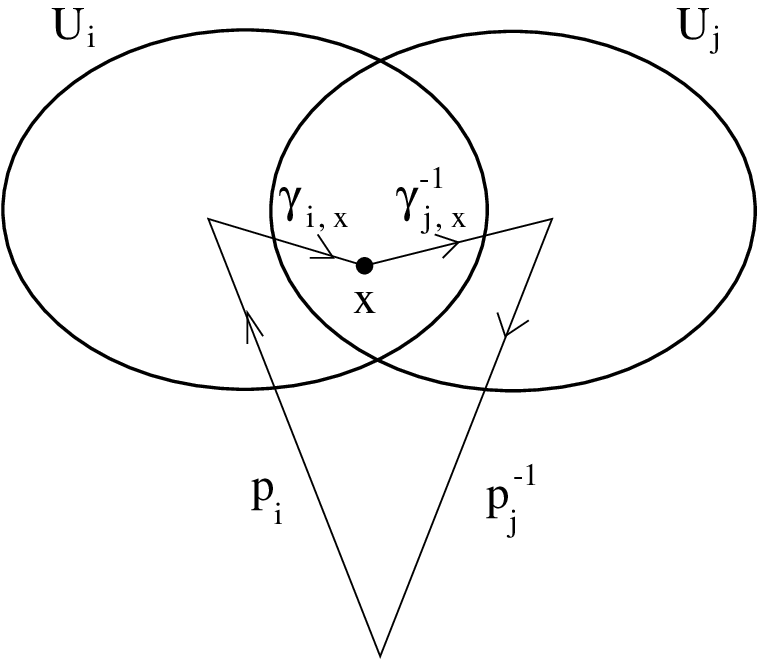}
}
\caption{transition function}
\label{transfun}
\end{figure}

Now we define the transition functions, $g_{ij}$, of the bundle 
$L$ corresponding to ${\cal H}$ by setting (see Fig.~\ref{transfun})
$$\ell_{ij}(x)=
p_{i,x}\star p_{j,x}^{-1}$$ 
and defining 
\begin{equation}
g_{ij}(x)= {\cal H}(\ell_{ij}(x))
\label{eq:transfun}
\end{equation}
\begin{Lem}
The transition functions $g_{ij}$ satisfy the cocycle condition.
\end{Lem}
{\bf Proof}: For $x$ in a triple overlap $U_{ijk}$, we have 
$g_{ij}g_{jk}g_{ki}(x)=1$,
since the product of the corresponding three loops is rank-$1$ homotopic to the
trivial constant loop, using the fact that we may cancel
products of the form $p\star p^{-1}$ up to rank-$1$ homotopy. Note that 
we have $g_{ji}=g_{ij}^{-1}$ by definition as well.\qed

\begin{figure}
\centerline{
\epsfbox{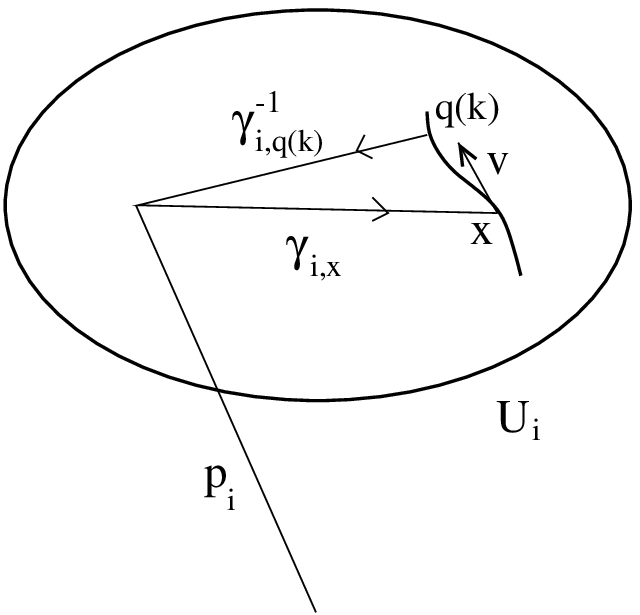}
}
\caption{connection}
\label{lconn}
\end{figure}

\begin{figure}
\centerline{
\epsfbox{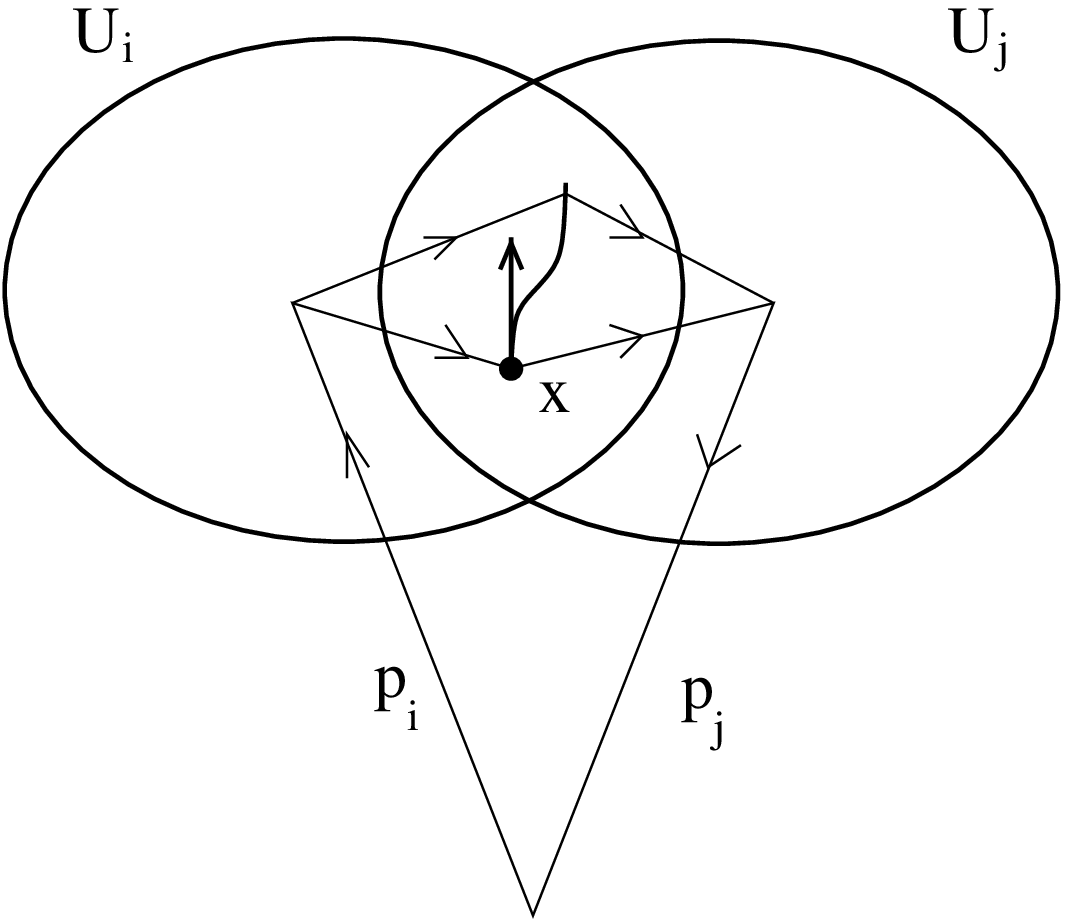}
}
\caption{Proof of Eq.~\ref{eq:lineconnproof}}
\label{lineconnpro}
\end{figure}

To define the $1$-forms $A_i$ corresponding to ${\cal H}$, 
let $x\in U_i$, $v\in T_xU_i$, and let 
$q: ]-\epsilon, \epsilon [ \rightarrow U_i$ be a smooth path, such that $q(0)=x$,
$\dot{q}(0)=v$.  Let $q_k$ denote a path following $q$ from $x$ to $q(k)$,
reparametrized at $t=0,1$ so as to belong to $P^\infty(M)$. Concretely we 
define $q_k(t)=q(\beta(t)k)$, where $\beta$ is the previously mentioned 
smoothly increasing 
function which is equal to $0$ on $[0,\frac{1}{3}]$ and equal to $1$ on 
$[\frac{2}{3},1]$. Note that $d/dk\, \left(q_k(1)\right)\vert_{k=0}=v$. 
Define the loop (see Fig.~\ref{lconn})
$$\ell_{i,q}(k)=p_{i,x}\star q_k \star p_{i,q(k)}^{-1}$$ 
and set 
$$f_{i,q}(k)={\cal
H}(\ell_{i,q}(k)).$$ 
Now we define:
$$
A_i(v)= -i\frac{d}{dk} \log f_{i,q}(k) \left. \right|_{k=0}.
$$
\begin{Lem}
The $1$-forms $A_i$ are well-defined and they satisfy 
$$i\left(A_j-A_i\right) = d\log g_{ij} $$ 
on double overlaps $U_{ij}$.
\end{Lem}
{\bf Proof}: From Fig.~\ref{lineconnpro} we have the equality
\begin{equation}
g_{ij}(q(k)) = f_{i,q}(k)^{-1}g_{ij}(x) f_{j,q}(k).
\label{eq:lineconnproof}
\end{equation}
Taking the derivative of the logarithm at $k=0$ we derive
$$
g_{ij}^{-1}\, d \, g_{ij}(v)= i(A_j(v)-A_i(v)).
$$
Now we introduce a new open set $U_j$ in the atlas and corresponding 
$\phi_j: U_j \rightarrow B(0,1)$, such that $U_j$ is centred around $x_j=x$ and contained 
in $U_i$. Such a pair $U_j$, $\phi_j$ may easily be constructed from $\phi_i$. Take the path from $*$ to the center of $U_j$ to be 
$p_j=p_{i,x}$. Now let $r$ be a second path in $U_i$ satisfying 
$r(0)=x$, $\dot{r}(0)=v$. Set 
\begin{eqnarray*}
A_i^{(q)}(v) & = &-i \frac{d}{dk} \log f_{i,q}(k) \left. \right|_{k=0} \\
A_i^{(r)}(v) & = &-i \frac{d}{dk} \log f_{i,r}(k) \left. \right|_{k=0} .
\end{eqnarray*}
Now
$$
A_j^{(q)}(v) - A_i^{(q)}(v) = A_j^{(r)}(v) - A_i^{(r)}(v),
$$
since $g_{ij}^{-1}\, d \, g_{ij}(v)$ is the evaluation of a $1$-form on a
vector and does not depend on which path is used. Let 
$\tilde{{\cal H}}: \pi_1^1(M, x) \rightarrow U(1)$ be the 
recentered holonomy at $x$ using the path $p_j$, i.e. 
$\tilde{{\cal H}}(\ell)= {\cal H}(p_j\star \ell  \star p_j^{-1})$ for any loop based at $x$.
Set $\tilde{f}_{j,q}(k)=q_k\star 
\gamma_{j,q(k)}^{-1}$. Then  
$$
A_j^{(q)}(v)= -i\frac{d}{dk} \log \tilde{{\cal H}}(\tilde{f}_{j,q}(k))\left. \right|_{k=0}=0
$$
where the last equality follows from Barrett's lemma to be proved in Sect.~\ref{BL}, since
$\tilde{f}_{j,q}(k)$ becomes the trivial loop at $x$ when $k$ goes to zero. Of 
course, by the same argument, $A_j^{(r)}(v)=0$, and thus $A_i^{(q)}(v)= 
A_i^{(r)}(v)$. 
\qed
The reconstruction described above involved choices of paths.

\begin{Lem} The reconstructed bundle and connection are independent of the
paths $p_i$ up to equivalence.
\end{Lem}
{\bf Proof}: Let $g'_{ij}$ and $A'_i$ be the reconstructed data using paths 
$p'_i$ instead of $p_i$. Define $h_i\in U(1)$ by $h_i={\cal H}(p_i\star (p'_i)^{-1})$, for each $i$.
Then 
$$g'_{ij}(x)=h_i^{-1} g_{ij}(x)h_j$$ 
\noindent and 
$$A'_i(v) = -i\frac{d}{dk} \log\, h_i^{-1}\star f_{i}(k)\star h_i
\left. \right|_{k=0}= A_i(v).$$ Thus the data $g'_{ij}$ and $A'_i$  
are equivalent to 
the original data $g_{ij}$ and $A_i$. \qed

The following result concerning the reconstructed connection will be used below.
\begin{Lem}
\label{calc}
Suppose $[a,b]\subset I$ is contained in $V_i$. Set 
$$\ell_{a,b}=p_{i,\ell(a)}\star 
\ell\left.\right|_{[a,b]}\star p_{i,\ell(b)}^{-1}.$$
Then  
$$
i \int_a^b \ell^\ast(A_i) = \log {\cal H} (\ell_{a,b}).
$$
\end{Lem}
{\bf Proof}: Set $\chi(k) = {\cal H} (\ell_{a,k})$. Now 
\begin{eqnarray*}
\frac{d}{dk} \log\chi(k) & = & \lim_{\epsilon\rightarrow 0} \log {\cal H} 
(\ell_{k,k+\epsilon})/\epsilon \\
& = & i A_i(\dot{\ell}(k)).
\end{eqnarray*}
Thus, integrating from $a$ to $b$, the result follows. \qed

It remains to show that the above assignments from bundle and connection 
data $L, {\cal A}$ to holonomies $\cal H$ and vice-versa are mutual inverses 
up to equivalence. Let ${\cal H}$ be the holonomy obtained from 
$(L,{\cal A})=(g_{ij},A_i)$.
The data reconstructed from this holonomy are given by:
\begin{eqnarray}
\tilde{g}_{ij}(x) & = & {\cal H}(\ell_{ij}(x))\label{eq:rec1}\\
i \tilde{A}_i(v) & = & \frac{d}{dk} \log {\cal H}
(\ell_{i,q}(k)) \left. \right|_{k=0}
\label{eq:rec2}
\end{eqnarray}
Now using the local formula for ${\cal H}$ in 
(\ref{eq:locforhol}) and the multiplicative property in (\ref{eq:mult}),  
and setting $h_i(x)=\H^{0i}(p_{i,x})$, we have
\begin{eqnarray*}
\tilde{g}_{ij}(x) & = & {\cal H}^{0i}(p_{i,x}) g_{ij}(x) 
{\cal H}^{j0}(p_{j,x}^{-1}) \\
& = & h_i(x)g_{ij}(x) h_j^{-1}(x)
\end{eqnarray*}
\noindent and
\begin{eqnarray*}
i \tilde{A}_i(v)  & = &  \frac{d}{dk} \log {\cal H}^{0i}
(p_{i,x}) 
{\cal H}^{ii} (q_k) 
{\cal H}^{i0}( p_{i,q(k)}^{-1}) \left. \right|_{k=0}\\
& = & \frac{d}{dk} i\int_I q_k^\ast (A_i) \left. \right|_{k=0}
+
\frac{d}{dk} \log {\cal H}^{i0}( p_{i,q(k)}^{-1}) \left. \right|_{k=0}\\
& = & i A_i(v) + d \log h_i(v).
\end{eqnarray*}
The final
equality follows from \begin{eqnarray*}
\frac{d}{dk} \int_0^1 q_k^\ast (A_i) \left. \right|_{k=0} & = & \frac{d}{dk} \int_0^k q^\ast (A_i) \left. \right|_{k=0}\\
& = & A_i(v).
\end{eqnarray*}
Thus $(\tilde{g}_{ij}, \,\tilde{A}_i)$ is equivalent to the original data 
$(g_{ij}, \, A_i)$. 

\begin{figure}
\vbox{\hskip1.5cm
\epsfbox{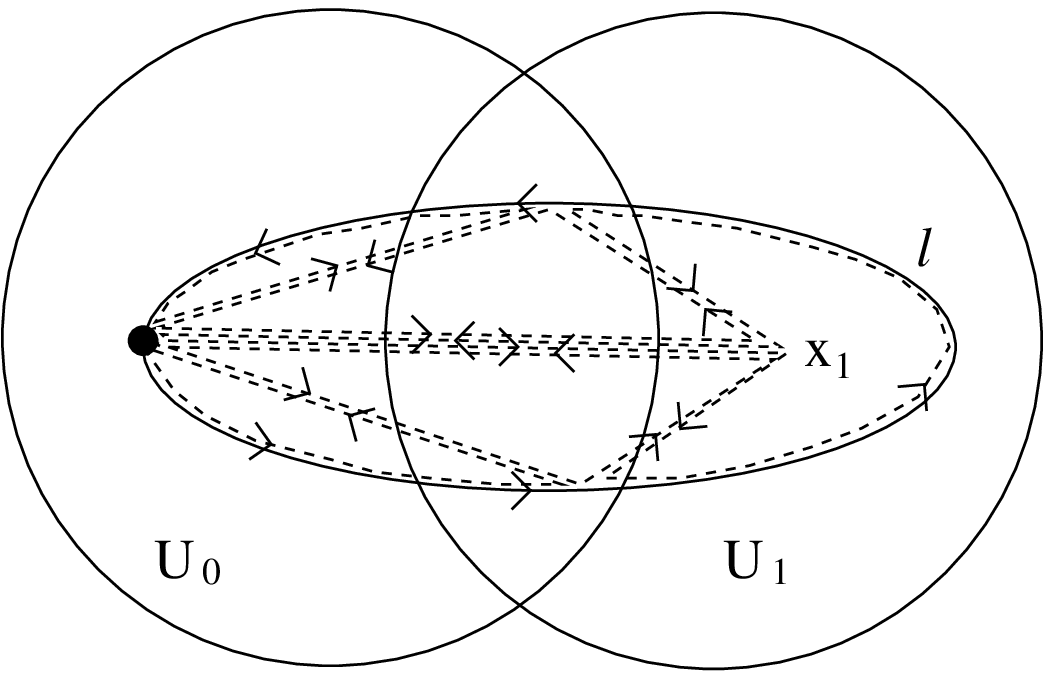}
}
\caption{$\H=\tilde\H$}
\label{biject}
\end{figure}

Conversely, let $(L,{\cal A})=(g_{ij},A_i)$ be the line bundle and 
connection obtained from the holonomy $\cal H$. Let
$\tilde{\cal H}$ be the holonomy obtained from these data. To show $\tilde{\cal
H}=\cal H$ we start with the local formula for $\tilde{\cal H}$:
\begin{equation}
\tilde{\cal H}(\ell) = 
\prod_{\alpha=1}^{N}\exp i\int_{I_{\alpha}}\ell^\ast(A_{\alpha})\cdot 
g_{\alpha,\alpha+1}(\ell(x_{\alpha})).
\label{Htilloc}
\end{equation}
Now using the definition of the transition functions (\ref{eq:rec1}) 
and Lem.~\ref{calc} 
all paths cancel except for the subpaths of $\ell$ between $x_{\alpha-1}$ 
and $x_{\alpha}$ (see Fig.~\ref{biject}), and therefore $\tilde{\cal H}(\ell) 
=
{\cal H}(\ell), \forall \ell$.

In conclusion we have the following 

\begin{Thm}
The above assignments define a bijection between bundles and
connections, modulo equivalence, on the one hand, and holonomies on the other.
\end{Thm}

\label{HLB}
\section{Thin higher homotopy groups}
Caetano and the second author~\cite{CP98} also defined higher (relative) smooth homotopy 
groups with homotopies of restricted 
rank, generalizing the rank-$1$ homotopy group $\pi_1^1(M)$ of the previous 
section. In the present article we will only need to consider the special case of homotopy groups 
relative to the base point, which simplifies the definition.

Let $M$ be a smooth manifold. Let $I^n$ denote the unit $n$-cube, with coordinates $t_i\in [0,1],\, 
i=1,\ldots, n$. 

\begin{Def} An $n$-loop is a smooth map $\gamma : I^n\rightarrow M$ such that, for some $0< 
\epsilon < 1/2$, 
$$\forall i=1, 
\ldots, n,\quad
t_i\in [0,\epsilon[\, \cup\, ]1-\epsilon, 1]\ \Rightarrow\ \gamma(t_1, 
\ldots , t_n) =\ast.$$
We denote the set of all $n$-loops by $\Omega_n^{\infty}(M)$.
\label{nloopdef}
\end{Def}

\begin{Rem} {\rm 
The above condition generalizes the ``sitting'' condition $\gamma(t)=\ast, \forall t\in 
[0,\epsilon[\cup ]1-\epsilon, 1] $ for loops in Sect.~\ref{HLB}. 
In~\cite{CP98} the weaker 
requirement
$$
(t_1\in [0,\epsilon[\cup ]1-\epsilon, 1]) \vee (t_i\in \{0,1\} 
\ \mbox{for some}\ i=2,\ldots,n)\Rightarrow\gamma(t_1, \ldots , t_n) =\ast
$$
was used.
Definition \ref{nloopdef} has the advantage that the $n$-loops can be multiplied smoothly (see below) 
along any of the $t_i$ directions, and not just along $t_1$.
}
\end{Rem}

\begin{Def} The product, $\gamma_1\star\gamma_2$, of two $n$-loops $\gamma_1$ and $\gamma_2$ is 
given by:
$$
\gamma_1\star\gamma_2(t_1,\ldots, t_n) =\left\{\begin{array}{ll}
\gamma_1(2t_1, t_2, \ldots, t_n), & t_1\in [0,1/2], \\
\gamma_2(2t_1-1, t_2, \ldots, t_n), & t_1\in ]1/2,1]. \\
\end{array}\right.
$$
The inverse, $\gamma^{-1}$, of an $n$-loop $\gamma$ is given by:
$$
\gamma^{-1}(t_1, \ldots, t_n) =\gamma(1-t_1,t_2, \ldots, t_n)
$$
\end{Def}

\begin{Rem} 
{\rm 
$\gamma_1\star\gamma_2$ is smooth because of the sitting condition in 
Def.~\ref{nloopdef}, which implies that $\gamma_1\star\gamma_2$ is constant in a neighbourhood of 
$t_1=1/2$.
}
\end{Rem}
We now define the thin homotopy relation, which was called intimacy relation 
by Caetano and Picken~\cite{CP98}. 
\begin{Def} Two $n$-loops $\gamma_1$ and $\gamma_2$ are said to be 
{\rm rank-$n$ homotopic} or {\rm thin homotopic}, denoted $\gamma_1 \stackrel{n}{\sim} \gamma_2$, if there exists $\epsilon >0$ and a 
homotopy 
$H: [0,1]\times I^n \rightarrow M$, such that:
\begin{itemize}
\item[1.] $t_i\in [0,\epsilon[\cup ]1-
\epsilon, 1]\ \Rightarrow\ H(s, t_1, \ldots ,t_n) =\ast,\quad i=1, \ldots, n$
\item[2.] $s\in [0,\epsilon[\ \Rightarrow\ H(s, t_1, \ldots ,t_n) = 
\gamma_1( t_1, \ldots ,t_n)$
\item[3.] $s\in  ]1-\epsilon, 1]\ \Rightarrow\ H(s, t_1, \ldots ,t_n) = 
\gamma_2( t_1, \ldots ,t_n)$
\item[4.] $H$ is smooth throughout its domain 
\item[5.] ${\rm rank} DH_{(s, t_1, \ldots ,t_n)} \leq n$ throughout its domain.
\end{itemize}
\label{rkndef}
\end{Def}
It is straightforward to show that $\stackrel{n}{\sim}$ is an equivalence relation. Let us denote the 
set of equivalence classes of $n$-loops in $M$ by $\pi_n^n(M,\ast)$, or just $\pi_n^n(M)$ when $\ast$ 
is understood.

\begin{Thm} $\pi_n^n(M,\ast)$ is an abelian group for $n\geq 2$.
\end{Thm}
{\bf Proof}: The product and inverse operations defined on $n$-loops descend to $\pi_n^n(M,\ast)$. 
The identity is the constant $n$-loop, which sends $I^n$ to $\ast$. The group properties are shown in 
the same way as for $\pi_n(M)$, with the modifications introduced by 
Caetano and Picken in~\cite{CP94,CP98} to accommodate smooth 
$n$-loops and homotopies. The group 
$\pi_n^n(M,\ast)$ is abelian by the standard geometric argument, since all homotopies involved are thin.
\qed

\begin{Rem} 
{\rm For $\dim M \leq n$, $\pi_n^n(M)=\pi_n(M)$ and there is nothing new. 
For $\dim M >n$ 
however, $\pi_n^n(M)$ is infinite-dimensional.
}
\end{Rem}

In the remainder of this paper we shall mainly be concerned with $\pi_2^2(M)$, the group of $2$-loops, or surfaces, modulo rank-$2$ homotopy.

\label{HH}
\section{Gerbe-holonomy}
Let $\G$ be a gerbe on $M$ given by a set of transition 
line-bundles 
$\Lambda_{ij}$ and trivializations 
$\theta_{ijk}$ , and 
let $\A$ be a gerbe-connection on $\G$ given by a set of connections, 
$\nabla_{ij}$, on $\Lambda_{ij}$ and local 
2-forms $F_{i}$ (see Sect.~\ref{LBG}).   
We first define the gerbe-holonomy of $(\G,\A)$, following 
Chatterjee~\cite{Ch98}. Let $s\colon I^2\to M$ be a 
2-loop, then the pull-back of $\G$, $s^*(\G)$, defines a gerbe on 
$I^2$. Since $I^2$ is two dimensional, the gerbe $s^*(\G)$ is trivial 
and we 
can choose an arbitrary trivialization, i.e., an object $\O$, with 
object connection. Let $\O$ be given by the line-bundles $L_i$, trivialized 
by the sections $\sigma_i$ over $V_i=s^{-1}(U_i)$, and let the 
object connection be given by 
$\nabla_i$. We can now define a global 2-form on $I^2$ by the formula 
$$\epsilon\vert_{V_i}=s^*(F_i)-\sigma^*_i(K(\nabla_i)).$$
Chatterjee~\cite{Ch98} calls $\epsilon$ an error 2-form of the object 
connection.
\begin{Def}\label{HGA} 
The {\rm holonomy} of $(\G,\A)$ around $s$, which we denote by $\H(s)$, is 
defined by 
$$\mbox{exp}\ i\int_{I^2}\epsilon.$$
\end{Def}
Chatterjee proves that $\H(s)$ is independent of all the choices we made.
\begin{Lem}{(\rm \cite{Ch98} Thm. 7.1.2)} The value of $\H(s)$ is 
independent of the choice of connection on the object, and of the choice of 
object.   
\end{Lem}
We now show that $\H$ is constant within each thin homotopy class. 
\begin{Lem} Let $s,s'\colon I^2\to M$ be two 2-loops. If $s$ and $s'$ 
are thin homotopic, then $\H(s)=\H(s')$.  
\end{Lem}
{\bf Proof}: Let $H\colon I^3\to M$ be a thin homotopy between $s$ and $s'$. 
Note that only the two faces corresponding to 
$s$ and $s'$ of $\partial H$ are 
mapped non-trivially to $M$, all other faces are mapped to the base-point. 
By the observation above and Stokes' theorem we get 
$$\H(s)\H(s')^{-1}={\rm exp}\ i\int_{\partial I^3}
\epsilon={\rm exp}\ i\int_{I^3}H^*(G).$$ 
Here $\epsilon$ is an error 2-form for an object of $H^*(\G)$ defined on 
$I^3$ and $G$ is the gerbe-curvature 3-form. Note that 
we can apply Stokes' theorem, to obtain the 
second equality, because 
$d\epsilon=H^*(G)$ (see~\cite{Ch98}). The last expression is equal 
to $1$ 
because the rank of the differential of $H$ is at most $2$ and $G$ is a 
3-form.\qed   
\noindent Since it is also clear from the definition that the holonomy of the 
product of two 2-loops equals the product of the holonomies around each one of them, 
we arrive at the following Lemma.
\begin{Lem} The gerbe-holonomy defines a group homomorphism
$$\H\colon\pi_2^2(M)\to U(1),$$ 
which only depends on $\G,\A$ up to equivalence.
\end{Lem}
{\bf Proof}: The first part of the claim is a corollary to the previous 
lemma. The second part follows directly from the definition of equivalence 
between gerbes with 
gerbe-connections by applying Stokes' theorem repeatedly.\qed

As we showed in Sect.~\ref{HLB} the holonomy of a line-bundle with connection 
is smooth in a precise sense. The same is true for gerbe-holonomies. We define 
a {\it smooth family of 2-loops} to be a map 
$\psi\colon U\subseteq\R^n\to\Omega^{\infty}_2(M)$ defined on an open subset 
$U\subseteq\R^n$ such that $\psi(x;t^1,t^2)=\psi(x)(t^1,t^2)$ is smooth 
on $U\times I^2$.
\begin{Def}\label{2hol} A {\rm 2-holonomy} is a group homomorphism $\H\colon\pi^2_2(M)\to U(1)$ such that for every smooth family of 2-loops $\psi\colon U\subseteq\R^n\to 
\Omega^{\infty}_2(M)$ the composite 
$$U\subseteq\R^n\stackrel{\psi}{\rightarrow}\Omega^{\infty}_2(M)\stackrel
{\rm proj}{\rightarrow}\pi_2^2(M)\stackrel{\H}{\rightarrow}U(1)$$
is smooth throughout $U\subseteq\R^n$.
\end{Def}
\begin{Lem} The gerbe-holonomy, $\H$, defines a 2-holonomy.
\end{Lem}
{\bf Proof}: This is an immediate consequence of the fact that $\H$ 
is defined by integration of a smooth 2-form.\qed
\begin{Lem}(Recentering a 2-holonomy) 
\begin{description}
\item{a)} Let $\ell\in\Omega^{\infty}(M,\ast)$ be a loop which is homotopic 
to the constant loop in $\ast$ via a homotopy $G\colon\ell\to\ast$. We 
define $P^2_2(M,\ell)$ as the group of all homotopies starting and ending at 
$\ell$ with fixed endpoints $\ast$ modulo thin homotopy. Given a group 
homomorphism $\tilde{\H}\colon P^2_2(M,\ell)\to U(1)$ we can {\rm recenter} 
$\tilde{\H}$ to obtain a 2-holonomy $\H\colon\pi_2^2(M)\to U(1)$ by defining 
$\H(s)=\tilde{\H}(G\star s\star G^{-1})$. We also say that we have {\rm 
recentered} $s$ by $G$.
\item{b)} Choose $m\in M$ and let $p\in P^{\infty}(M)$ be a path from 
$\ast$ to $m$. Given a 2-holonomy $\H\colon\pi_2^2(M,\ast)\to U(1)$ we 
can {\rm recenter} $\H$ to obtain a 2-holonomy $\tilde{\H}\colon\pi_2^2(M,m)
\to U(1)$ in the following way: choose $s\in\Omega^{\infty}_2(M,m)$ and let 
$s'$ be the 2-path $\mbox{Id}_p\star s\star\mbox{Id}_{p^{-1}}$. Now 
recenter $s'$ by the thin homotopy between the constant loop at $\ast$ and 
$p\star Id_{m}\star p^{-1}$, denoted by $G_{p}$. Note that for this 
recentering we have to use composition of 2-paths via the second coordinate. 
Thus we have obtained a 2-loop $\tilde{s}\in\Omega^{\infty}_2(M,\ast)$ and 
therefore we can define $\tilde{\H}(s)=\H(\tilde{s})$. We also say that 
we have {\rm recentered} $s$ by $p$.
\end{description}
\end{Lem}
{\bf Proof}: Both in a) and b) one only has to check that recentering is 
well defined modulo thin homotopy, which is straightforward.\qed

Before going on to the next section it is worthwhile to have a look at a more 
concrete formula for $\H$. Analogously to what we did for 
line-bundles with connections in Sect.~\ref{HLB}, we can define 
$\H$ completely in terms of the \v{C}ech cocycle 
$g_{ijk}$, the 
0-connection $A_{ij}$ and the 1-connection $F_i$. 
Let $s\colon I^2\to M$ 
be a 2-loop in $M$ as before. Let $\V=\left\{V_i\ \vert\ i\in J\right\}$ be the covering of $I^2$ obtained by taking the inverse image via $s$ of all open 
sets in the cover $\U$. We define a {\it grid} on $I^2$ 
to be a rectangular subdivision of $I^2$. A rectangle in the grid 
is denoted by $R_{i}$, an edge by $E_{ij}$ and a vertex by $V_{ijkl}$. The 
edges are all oriented from left to right and from bottom to top, the 
rectangles have the counter-clockwise orientation. Now take the grid 
sufficiently fine so that each rectangle $F_i$ lies in at least one 
open set $V_i$. Note that this is possible because 
each covering of a compact metric space has a Lebesgue number. Note also that 
in general a small rectangle can be contained in more than one open set. A 
particular choice of open set for each rectangle we call a {\it labelling} of 
the grid. We always assume that the labels of the rectangles at the boundary 
are equal to $0$, and we fix a choice of $V_0=s^{-1}(U_0)$, such that 
$U_0$ contains 
the base point of $M$. Let 
$\epsilon\in\Omega^2(I^2)$ be the error 2-form of an object for $s^*(\G,\A)$, 
given by a family of line-bundles $L_i$, with 0-connection given by a family 
of 1-forms $A_i\in\Omega^1(V_i)$. We write 
$\epsilon\vert_{V_i}=s^*F_i-dA_i$. We recall the identities
$$i(A_j-A_i)=is^*A_{ij}+d\log\lambda_{ij},$$ 
where $\lambda_{ij}\colon V_{ij}\to U(1)$ satisfies  
$$\lambda_{ij}\lambda_{jk}\lambda_{ki}=s^*g_{ijk}.$$
Now pick a labelling of the grid. 
According to Def.~\ref{HGA} we have 
\begin{equation}
\begin{array}{lll}
\H(s)&=&\exp\ i\int_{I^2}\epsilon\\[15pt]
&=&{\displaystyle\prod_{\alpha}}\exp\ i\int_{R_{\alpha}}
(s^*F_{\alpha}-dA_{\alpha})\\[15pt]
&=&{\displaystyle\prod_{\alpha}}\exp\ i\int_{R_{\alpha}}s^*F_{\alpha}\cdot
{\displaystyle\prod_{\alpha,\beta}}\exp\ i\int_{E_{\alpha\beta}}s^*
A_{\alpha\beta}\\[15pt]
&&\times\,{\displaystyle\prod_{\alpha,\beta,\gamma,\delta}}g_{\alpha\beta\gamma}(s(V_{\alpha\beta\gamma\delta}))g_{\alpha\delta\gamma}^{-1}(s(V_{\alpha\beta\gamma\delta})).
\label{eq:gholfor}
\end{array}
\end{equation} 

\begin{figure}
\centerline{
\epsfbox{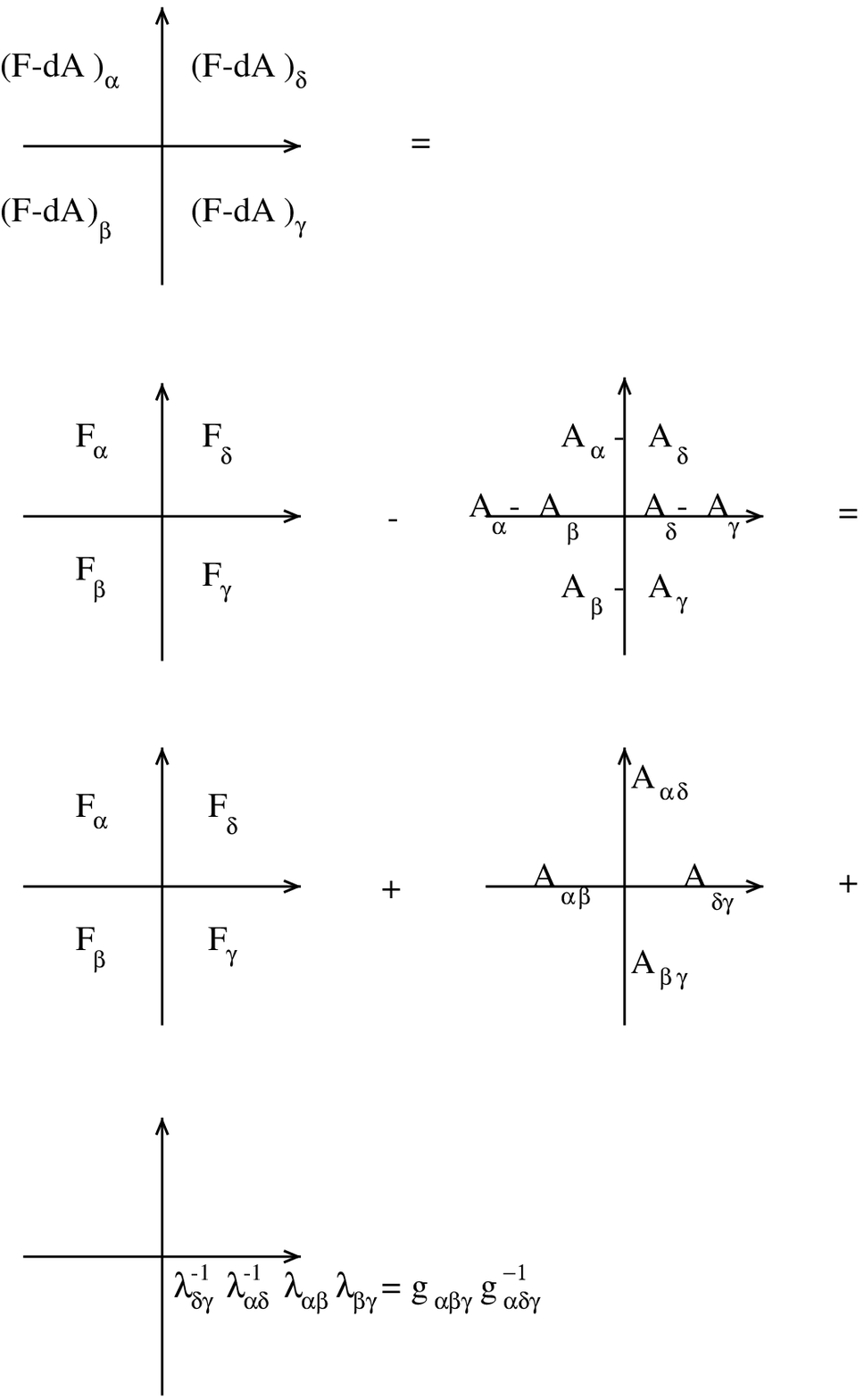}
}
\caption{concrete formula for gerbe-holonomy}
\label{concrete}
\end{figure}
\noindent The last two products are to be taken over the labels of 
contiguous faces in the grid only and in such a way that 
each face, edge and vertex appears only once. The convention for the 
order of the labels is indicated in Fig.~\ref{concrete}. 
Formula (\ref{eq:gholfor}) follows from applying Stokes' theorem 
repeatedly. 
Note that Chatterjee's results for the global definition of the 
gerbe-holonomy~\cite{Ch98} and the equalities above show that the value of 
$\H(s)$ is independent of the choice of grid and its labelling.

We have written down the explicit formula for gerbe-holonomy using a grid. 
Of course it is possible to obtain an analogous formula using other 
subdivisions of the unit square, for example, triangulations of $I^2$. The 
idea is exactly the same, but one has to take into account the different 
valencies of the vertices. 

In Sect.~\ref{Rec} we show that every 2-holonomy is the gerbe-holonomy of 
some gerbe with gerbe-connection.

\label{HG}
\section{Parallel transport in gerbes}
Let us now explain parallel transport for Abelian 
gerbe-connections in general along arbitrary homotopies between arbitrary 
loops. For this discussion we 
assume some basic knowledge about groupoids and 2-groupoids, which can be 
obtained by reading~\cite{Ba972,BD95,Mack87} for example. We should warn 
the reader that 
this section is different in flavour from the rest of the paper. Whereas 
we have followed a down-to-earth approach in the rest of the paper, here 
some readers might feel that we are trying to make up for that by being 
unnecessarily sophisticated. We have two arguments in our defense. First of 
all we do not know any simpler way of formulating the parallel transport 
of gerbes satisfactorily. Secondly, we believe that category theory is already 
at the heart of gerbes, since the original definition of gerbes is 
in terms of sheaves of categories~\cite{Bry93}. Being at the heart of gerbes, 
we ought to understand the category theory that is involved a little better. 
Our formulation of the 
parallel transport of a gerbe on $M$ shows the relation with what can 
be called the {\it thin homotopy $2$-type} of $M$. Following ideas 
expressed by Grothendieck, homotopy theorists and category theorists 
have been endeavouring to define {\it weak $n$-categories}, of which a 
special sub-class, the {\it weak $n$-groupoids}, should model homotopy 
$n$-types 
of topological spaces. For an overview of $n$-category theory 
see~\cite{Ba972,BD95}. In~\cite{BD95} Baez and Dolan sketch the possible 
relevance of $n$-categories for the formulation of 
{\it Topological Quantum Field Theories} (TQFT's). 
Following these ideas the first author of this article 
investigated the possibility of defining four-dimensional TQFT's using 
monoidal 2-categories~\cite{Ma99,Ma00}. Our formulation of parallel transport 
of gerbes in this section is also a first attempt to see 
if there is any link between 
$n$-categories and TQFT's on the one hand and differential geometry on the 
other hand. 
If only for 
this reason, we already feel that the effort of penetrating the relatively 
unfamiliar language of monoidal categories and 2-categories is not wasted. 
Due to this higher level of sophistication we cite results from the literature 
without giving direct proofs here, so that the full emphasis is placed on 
the change of language and not on mathematical detail. For mathematical 
detail we refer to Sect.~\ref{Rec}, where we work out the simply-connected 
case in the more familiar language of groups and group homomorphisms. At 
the end of this section we indicate how an explicit proof in the general 
case can be deduced from the results in Sect.~\ref{Rec}. 

It is enlightening to go back to the case of principal bundles with 
connections first.  
Let $P\stackrel{p}{\to} M$ be a principal $G$-bundle on a 
not-necessarily-connected manifold $M$ and let $\omega$ be a connection in 
$P$. 
Suppose in that case that we are trying to 
describe the parallel transport of $\omega$. It is well known in 
homotopy theory that 
for general manifolds it is best to work with the {\it path 
groupoid} rather than the fundamental group, because the latter requires 
the choice of a basepoint in a certain connected component. Therefore it is 
natural to employ the thin path groupoid, $P_1^1(M)$ (see Sect.~\ref{HLB}) 
to describe 
parallel transport. To formulate parallel transport along paths in terms of 
Lie groupoids and functors we also have to associate a groupoid to the 
bundle $P$. This is 
a well known construction which goes back to Ehresmann's work, and details can 
be found in Mackenzie's book~\cite{Mack87}, for example. 

Let us sketch the 
construction of this groupoid, which we denote by $G(P,M)$. 
\begin{Def}
The objects of $G(P,M)$ are simply the points in $M$. The 
set of all morphisms is given by the manifold $P\times P/G$, where $G$ 
acts by $(x,y)s=(xs,ys)$. We denote the equivalence classes by $[x,y]$, 
where we consider this to be a morphism from $p(x)$ to $p(y)$, the 
opposite of Mackenzie's~\cite{Mack87} convention. The 
composite $[x,y][w,z]$ is defined by $[xs,z]$, where $w=ys$ in the fibre 
over $p(y)=p(w)$. The identity morphism for $m\in M$ is of course $[x,x]$, 
where $x\in p^{-1}(m)$ is arbitrary. Finally the inverse of $[x,y]$ is 
given by $[y,x]$.
\end{Def}
It is easy to check that this defines a groupoid, and 
it can be shown~\cite{Mack87} that it is a locally trivial Lie groupoid.  
This last statement means that the operations just defined are smooth and 
that for every point in $M$ there is an open neighborhood, say $U$, such that 
the restriction of this groupoid to $U$ is isomorphic to 
$U\times G\times U$, with the trivial groupoid structure. Conversely one can 
construct a principal bundle from a locally trivial Lie groupoid by taking 
as the total space all morphisms with a fixed, but arbitrary, source, 
and as the projection the target map. These constructions are each other's 
inverses up to isomorphism.      

The connection $\omega$ now gives rise to a functor ${\cal P}\colon 
P_1^1(M)\to G(P,M)$ which is the identity on objects, i.e. points of $M$, 
and is smooth in the 
sense of Sect.~\ref{HLB}. 
\begin{Def}
The {\it P(arallel) T(ransport) functor}, denoted ${\cal P}$, of $\omega$ is 
defined on objects by ${\cal P}(m)=m$, for all $m\in M$. For a given 
path $q$ in $M$, we define ${\cal P}(q)=[x,y]$, 
where $x\in p^{-1}(q(0))$ is arbitrary and $y\in p^{-1}(q(1))$ is obtained 
from $x$ by parallel transport along the path $q$. Since parallel transport 
only depends on the thin homotopy class of $q$, the functor ${\cal P}$ 
is well-defined.
\end{Def}
Conversely, any such functor yields a 
{\it path-connection} in $G(P,M)$ 
in the sense of Def.~7.1 in~\cite{Mack87}, and is therefore equivalent to a 
connection 
in $P$. The path-connection for a given PT-functor, $\cal P$, is easy to 
describe: for any path $q$ in $M$, the path-connection yields the 
path in $G(P,M)$ 
that is given by $k\mapsto {\cal P}(q_k)$, where $q_k$ is as 
in Sect.~\ref{HLB}. 
The definition of ${\cal P}$ and the observations above cannot be found in 
the literature on Lie groupoids, but they are an 
immediate consequence of Barrett's construction~\cite{Ba91}, so we will not 
spell them out here. It is easy to see that two connections in $P$ are 
gauge-equivalent precisely if the corresponding $PT$-functors are naturally 
isomorphic. Note that the reconstruction result is less powerful than 
Barrett's original 
theorem, because the bundle is already part of the PT-functor in the form of its target. However, the upshot is that we can deal 
with the non-connected case directly by not fixing the choice of a 
base-point. 

Back to gerbes again. Brylinski~\cite{Bry93} (see Hitchin~\cite{Hi99} too), 
explains how 
a gerbe on $M$ with a gerbe-connection gives rise to a line-bundle on 
$\Omega(M)$ with an ordinary connection. This line-bundle has the special 
property that it is ``multiplicative'' with respect to the composition of 
loops (Prop. 6.2.5 in~\cite{Bry93}). The connection in this line-bundle has the special 
property that it yields a parallel transport over ``cylinders'' which 
does not depend on the way the cylinder is made up out of a path of loops 
but only on the surface itself. This leads to the idea that a gerbe, $\cal G$, 
with gerbe-connection, $\cal A$, on $M$ yields a groupoid on the 
{\it thin loop groupoid}, $L^1_1(M)$, which is the subgroupoid of $P^1_1(M)$ 
with only loops. If $M$ is path-connected (something which we will 
assume for the rest of this section), one can fix a basepoint, $*\in M$, 
and work 
over $\pi_1^1(M,*)$. Before 
explaining the parallel transport functor for gerbes, let us first have a 
closer look at the construction of this line-bundle with connection over 
$\pi_1^1(M,*)$.

We follow Ch.~6 in Brylinski's book~\cite{Bry93} mostly, but we consider 
$\Omega(M,*)$ to be smooth using smooth families of loops rather than trying 
to 
define an infinite-dimensional manifold structure on it. Let $\cal G$ be a 
gerbe on 
$M$ and $\cal A$ a gerbe-connection. We first construct the line-bundle on 
$\Omega(M,*)$, as Brylinski does, and then show that it projects to a 
line-bundle on $\pi_1^1(M,*)$. Both bundles we denote by $L^{\G}$. 
\begin{Def}
The total 
space of $L^{\G}$ is given by the set of equivalence classes of quadruples 
$(\gamma,F,\nabla,z)$, where $\gamma\in\Omega(M,*)$, $F$ is an object for 
$\gamma^*\G$ on the circle $S^1$ with $\nabla$ an object connection 
in $F$, and finally $z\in\Co^*$. Note that we use the word ``object'' 
as defined in Sect.~\ref{LBG}, so it really is given by a set of data, as is $\nabla$.  
The equivalence relation is generated 
by
\begin{enumerate}
\item $(\gamma,F,\nabla,z)\sim (\gamma,F',\nabla',z)$ if $(F,\nabla)$ and 
$(F',\nabla')$ are isomorphic pairs of objects with object-connections. 
\item $\left(\gamma,F,\nabla+\alpha,z\right)\sim 
\left(\gamma,F,\nabla,z\exp\left(-\int_0^1\alpha\right)\right)$ for any 
complex valued 1-form $\alpha$ on $S^1$.
\end{enumerate}   
The action of $\Co^*$ on $L^{\G}$ is given by 
$(\gamma,F,\nabla,z)w=(\gamma,F,\nabla,zw)$. Brylinski (Prop.6.2.1. 
in~\cite{Bry93}) proves that there is a unique smooth structure on 
$L^{\G}$ such that 
\begin{enumerate}
\item The projection $(\gamma,F,\nabla,z)\mapsto\gamma$ defines a smooth 
principal $\Co^*$-bundle.
\item For any $\gamma\in\Omega(M,*)$, any open contractible neighborhood of 
the origin $U\subset\R^n$, and any smooth family of loops $\Gamma\colon 
U\to\Omega(M,*)$ such that $\Gamma(0)=\gamma$, let $(F,\nabla)$ be an 
object with object-connection on $\Gamma(U)\subset M$. The map  
$$\sigma(\Gamma(x))=(\Gamma(x),\Gamma(x)^*F,\Gamma(x)^*\nabla,1)$$
defines a smooth local section of $L^{\G}$.
\end{enumerate}
\end{Def}
Recall that two objects with, in this case necessarily flat, 
object-connections of 
$\gamma^*\G$ on $S^1$ always differ by a line-bundle with a flat connection on 
$S^1$. The two objects with connections are isomorphic if and only if 
the flat connection in the line-bundle has trivial holonomy (in which case the 
line-bundle is necessarily trivializable as well). The fibre over 
a loop is now acted upon by line-bundles on $S^1$ with flat connections, 
whose isomorphism classes are the elements of 
$\mbox{Hom}\left(\pi_1(S^1)=\Z,U(1)\right)\cong U(1)$. This is 
Hitchin's~\cite{Hi99} description of $L^{\G}$. As remarked by Brylinski 
(Prop.~6.2.5. \cite{Bry93}) the fibres of $L^{\G}$ have a 
multiplicative property. 
Given the points $(\gamma,\gamma^*F,\gamma^*\nabla,1)$ in $p^{-1}(\gamma)$ 
and $(\mu,\mu^*F,\mu^*\nabla,1)$ in $p^{-1}(\mu)$, the 
product becomes 
$(\gamma\star\mu,(\gamma\star\mu)^*F,(\gamma\star\mu)^*\nabla,1)$ in 
$p^{-1}(\gamma\star\mu)$. Here we have taken an object $F$ and an 
object-connection $\nabla$ which are defined on the image of $\gamma\star\mu$.

The connection on $L^{\G}$ can now be defined. Let $v$ be any tangent 
vector to $\Omega(M,*)$ at $\gamma$. Let $\sigma$ be the section 
defined above and let $\epsilon$ be the error 2-form (see Sect.~\ref{HG}) of 
$(\gamma^*F,\gamma^*\nabla)$ with the same notation as above. 
\begin{Def}
The covariant derivative in $L^{\G}$ is defined by 
$$\frac{D_v\sigma}{\sigma}=i\int_0^1\epsilon_{\gamma(t)}\left(\dot{\gamma}(t),
v(\gamma(t))\right)dt.$$
\end{Def}
The curvature of this connection is equal to 
$$K_{\gamma}(D)(u,v)=\int_0^1\Omega_{\gamma(t)}\left(\dot{\gamma}(t),
u(\gamma(t)),v(\gamma(t))\right)dt,$$
where $\Omega$ is the gerbe-curvature 3-form on $M$. It is also easy to 
describe the parallel transport of the connection $D$ 
along a ``cylinder''. Let $H\colon 
I^2\to M$ be a 
homotopy between two loops $\gamma$ and $\mu$. Choose an object $F$ and an 
object connection $\nabla$ for $H^*(\G)$ on $I^2$ (or $S^1\times I$). Let $\epsilon$ 
be the error 2-form for this object-connection. 
\begin{Def}
The parallel transport along $H$ is given by 
$${\cal P}(H)\left(\gamma,\gamma^*F,\gamma^*\nabla,1\right)=
\left(\mu,\mu^*F,\mu^*\nabla,1\right)\exp\left(\int_{I^2}\epsilon\right).$$ 
\end{Def}
Two observations show that $\left(L^{\G},D\right)$ 
projects to a bundle on $\pi_1^1(M,*)$. Firstly the formula for 
parallel transport clearly shows that it is compatible with the 
multiplication. Secondly, as $\epsilon$ is a 2-form, the parallel transport 
along a thin homotopy is trivial, so there is a unique way to 
identify the fibres of $L^{\G}$ which lie over loops that are thin homotopic.
By abuse of notation we denote this line-bundle over $\pi_1^1(M,*)$ with connection also 
by $(L^{\G},D)$. Furthermore, Brylinski (\cite{Bry93}, Thm.~6.2.4.(3)) 
shows that whenever two homotopies $G,H\colon I^2\to M$ between a given 
pair of paths are homotopic themselves by a homotopy $J\colon I^3\to M$, then 
the parallel transport around $GH^{-1}$ is given by 
$$\int_{I^3}J^*\Omega.$$ 
Thus we see that the parallel transport along $G$ equals the parallel 
transport along $H$ if they are thin homotopic, because in that case $J$ can 
be chosen to have rank $\leq 2$ everywhere.

One can now define the groupoid 
$L^{\G}\times L^{\G}/U(1)$ over $\pi_1^1(M,*)$ as we showed above, and 
this groupoid is equipped with 
a {\em monoidal structure} in the following sense. The equivalence 
classes of the loops can be multiplied, or tensored, 
$[\gamma]\star[\mu]=[\gamma\star\mu]$. Given $\alpha,\beta,
\gamma,\mu\in\Omega(M,*)$, and given two morphisms in $L^{\G}\times 
L^{\G}/U(1)$, say $[a,b]\colon[\alpha]\to[\beta]$ and 
$[c,d]\colon[\gamma]\to[\mu]$, then one can 
tensor the morphisms to get 
$[a,b]\star[c,d]=[ac,bd]\colon[\alpha\star\gamma]\to[\beta\star\mu]$, 
where $ac$ is defined 
by the product of the fibres as described above. This multiplication of the 
morphisms is compatible with the composition of morphisms in the following 
sense: for any quadruple of morphisms we have the equality 
$$([a,b]\star[c,d])([e,f]\star[g,h])=([a,b][e,f])\star([c,d][g,h]),$$ 
whenever both sides of the equation are defined. This equation, which is an 
example of the {\it interchange law} for monoidal categories, is easily 
checked and only holds because $U(1)$ is Abelian. With respect to this monoidal 
structure the objects and the morphisms have inverses and there is a unit 
object, $[c_*]$, where $c_*$ is the constant loop at $*$, and a 
unit morphism, $[a,a]$ for any $a$ such that 
$p(a)=[c_*]$. Altogether we propose to call this structure a 
{\it Lie 2-group}, 
which is justified by the fact that all operations involved are smooth in 
the appropriate sense and that it can be seen as a Lie 2-groupoid with 
only one object. (This is analogous to the statement that a Lie groupoid 
with one object is nothing but a 
Lie group. For this kind of general remark about $n$-categories 
see~\cite{Ba972,BD95,BD98} for example.) Without the smoothness condition this kind of 
groupoid goes under a variety of names in the homotopy literature. 
Yetter~\cite{Ye93} calls them {\it categorical groups}, for example. 
\begin{Def}
Let $G(\G,\A,M)$ be the Lie $2$-group given by 
$L^{\cal G}\times L^{\cal G}/U(1)$ over $\pi_1^1(M)$, as defined above.
\end{Def} 
Let us now define the {\it thin Lie 2-group of cylinders}, denoted by 
$C_2^2(M,*)$. We define $C_2^2(M,*)$ as the quotient of a non-strict monoidal 
groupoid $C(M,*)$ by a {\it normal monoidal subgroupoid} $N(M,*)$. 
\begin{Def} The 
objects of 
$C(M,*)$ are the elements of $\Omega(M,*)$. The morphisms are thin homotopy classes 
of homotopies between loops, through based loops. 
It is clear that this forms a groupoid under the 
obvious composition of homotopies. There is also a monoidal structure on 
$C(M,*)$ defined by the composition of loops and the corresponding composition 
of homotopies. For clarity, we refer to the monoidal composition of 
homotopies as {\it horizontal composition}, and 
write $\star$, and the other we call 
{\it vertical} and indicate by simple concatenation. We denote the vertical 
inverse of a homotopy $H$ by $H^{-1}$ and the horizontal inverse by 
${H}^{\leftarrow}$. 
\end{Def}
Recall that for 
homotopies of the trivial loop to itself 
both compositions are the same up to thin homotopy, which is why 
$\pi_2^2(M)$ is abelian. Concretely this follows from the interchange 
law which states that $(G_1\star H_1)(G_2\star H_2)$ is thin homotopic to 
$(G_1G_2)\star(H_1H_2)$, whenever both composites can be defined. 
This way the groupoid $C(M,*)$ is a {\it weak} monoidal groupoid, 
with weak inverses for the objects because loops only form a group up to 
thin homotopy. Instead of using the abstract {\it strictification 
theorem}~\cite{MP85}, which is not very practical for the concrete application 
to gerbe-holonomy, we ``strictify'' $C(M,*)$ by hand by dividing out by the 
monoidal subgroupoid of only the thin homotopies, $N(M,*)$. Dividing 
out by a monoidal subgroupoid is only well-defined if the 
following conditions are satisfied, in which case we call it {\it normal}:
\begin{enumerate}
\item For any $\gamma\in\Omega(M,*)\colon\quad 1_{\gamma}\in N(M,*)$.
\item For any $\gamma,\mu\in\Omega(M,*)$, any homotopy 
$G\colon\gamma\to\mu$ and any {\bf thin} homotopy $H\colon\mu\to\mu$, the 
{\bf thin} homotopy class of $GHG^{-1}\colon\gamma\to\gamma$ belongs to 
$N(M,*)$.
\item For any $\gamma,\mu\in\Omega(M,*)$, any homotopy $G\colon\gamma\to\gamma$ and any 
{\bf thin} homotopy $H\colon\mu\to\mu$, the {\bf thin} homotopy class of 
$G\star H\star G^{\leftarrow}\colon\gamma\star\mu\star\gamma^{-1}\to 
\gamma\star\mu\star\gamma^{-1}$ belongs to $N(M,*)$.  
\end{enumerate}
It is easy to check that the conditions above are the right ones for 
our construction of the quotient monoidal groupoid, which we 
explain below, to be well-defined. The only reference for the definition of 
a normal monoidal subgroupoid that we know of is~\cite{BW99} (who give the 
more general definition of a normal monoidal subcategory), but it might be 
that it can be found in earlier papers on monoidal categories and 
2-categories. We suspect that this definition goes back to the time when 
monoidal categories were defined for the first time~\cite{Be63}, but we 
have been unable to find a precise written reference in the older literature. 

\begin{Lem}
\label{wd}
The monoidal subgroupoid $N(M,*)$ is normal in $C(M,*)$.
\end{Lem}
{\bf Proof} The first condition is obviously satisfied. 

We prove that the 
second condition holds. Denote the group of homotopies 
$\gamma\to\gamma$ in $C(M,*)$ by $C(M,*)(\gamma)$. This group is isomorphic to 
$C(M,*)(\gamma\star\mu)$, for any 
$\mu\in\Omega(M,*)$ (this is well known, 
see for example~\cite{HKK01}). The isomorphism, which clearly preserves 
thinness, 
is given by 
$$\phi\colon G\mapsto G\star 1_{\mu}.$$ Under $\phi$ the homotopy $GHG^{-1}\colon\gamma\to\gamma$ 
is mapped to $(G\star 1_{\mu})(H\star 1_{\mu})(G^{-1}\star 1_{\mu})
\colon\gamma\star\mu\to\gamma\star\mu$. Clearly this is thin homotopic to 
$(G\star 1_{\mu})(1_{\mu}\star H)(G^{-1}\star 1_{\mu})$. By the interchange 
law the latter is thin homotopic to 
$1_{\gamma}\star H\colon\gamma\star\mu\to\gamma\star\mu$ 
which is clearly a thin homotopy whenever $H$ is thin.

The third condition can be proved by a similar argument.\qed
We can now define the quotient groupoid $C_2^2(M,*)=C(M,*)/N(M,*)$. 
\begin{Def}
The objects 
of $C_2^2(M,*)$ are the elements of $\pi_1^1(M,*)$, which we temporarily 
denote by $[\gamma]$. 

For any $\alpha,\beta,\gamma,\mu\in\Omega(M,*)$ and 
for any $G\in C(M,*)(\alpha,\beta)$ and $H\in C(M,*)(\gamma,\mu)$, we say that 
$G$ and $H$ are equivalent if there exist thin homotopies $A\in 
N(M,*)(\alpha,\gamma)$ and $B\in N(M,*)(\beta,\mu)$ such that 
$AHB^{-1}\stackrel{2}{\sim} G$. 

The morphisms between 
$[\gamma]$ and $[\mu]$ are 
the equivalence classes of 
$$\bigcup_{\alpha,\beta}\left\{C(M,*)(\alpha,\beta)\ \vert\ 
[\alpha]=[\gamma],\ [\beta]=[\mu]\right\}$$ modulo this equivalence relation.

The composition and the monoidal structure descend to the quotient precisely 
because $N(M,*)$ is normal. The smooth structure is defined by smooth families 
of loops and smooth families of homotopies. 
\end{Def}
In this way we 
have obtained a second example of a Lie 2-group.

Before going on we should explain a small technical gap that we have not been 
able to bridge yet. Although Lem.~\ref{wd} is strong enough to conclude that 
$C_2^2(M,*)$ is well-defined, it is too weak to prove that $C(M,*)$ and $C_2^2(M,*)$ 
are equivalent as Lie $2$-groups. To establish that equivalence one would have 
to prove that, for any $\gamma\in\Omega(M,*)$, the group $N(M,*)(\gamma)$ is 
trivial. We conjecture that this true, but have not been able to prove it 
for any manifold other than $\R^n$, where it is almost immediate. As it 
stands we 
cannot prove that $C_2^2(M,*)$ is a true strictification of $C(M,*)$, 
but it is clear that $C_2^2(M,*)$ suits our purpose in this paper 
very well indeed. 
In~\cite{HKK00} the reader can find a more restricted notion of thin 
homotopy for 
which the conjecture can be shown. However, this notion of thin homotopy 
does not seem to be suited to the smooth context of 
parallel transport. Besides being right 
for the formulation of parallel transport, our notion of 
thinness has the extra benefit that it can be generalized immediately to 
homotopies of any dimension, which is necessary if one wants to formulate and 
prove the analogues of our results for $n$-gerbes in general.

From everything above it now follows that the following theorem holds:
\begin{Thm}\label{2group}The pair $(\G,\A)$ gives rise to a 
smooth PT-functor of Lie $2$-groups 
$${\cal P}\colon C_2^2(M,*)\to G(\G,\A,M),$$
over the identity on $\pi_1^1(M,*)$.

Conversely, given any line-bundle $L$ on $\pi_1^1(M,*)$ with the 
multiplicative property as described above, we obtain a Lie 2-group, $G(L,M)$.
Any smooth functor ${\cal P}\colon C_2^2(M,*)\to G(L,M)$ over the 
identity on $\pi_1^1(M,*)$ gives rise to a 
gerbe-connection, unique up to gauge-equivalence, in the gerbe associated to 
$G(L,M)$, such that $\cal P$ is the PT-functor for that connection.
\end{Thm}
One can wonder about the meaning of the second part of this theorem. Where 
would one get such line-bundles on the loop space? Our construction of 
the line-bundle on $\pi_1^1(M,*)$ only depends on the gerbe and the 
$0$-connection. Thus in this formulation the $1$-connection is really 
the new information contained in the PT-functor. It is interesting to 
compare this to the case of line-bundles described in the beginning of this 
section, where the line-bundle is already given in the form of a Lie groupoid 
but the connection is determined by the PT-functor.

An explicit proof of the second part of Thm.~\ref{2group} in which one 
recovers the \v{C}ech 2-cocycle of the gerbe and the local 1- and 2-forms of 
the 
gerbe-connection from the parallel transport is easily deduced from our 
results in Sect.~\ref{Rec}, where we work out the simply-connected case 
explicitly and in great detail. There we close up all $2$-paths to obtain 
$2$-loops based at $\ast$ using certain auxiliary homotopies called 
$P_{ij}$, which exist when $M$ is $1$-connected. However, one can leave 
them out in general to obtain 
a direct proof of the second part of Thm.~\ref{2group}. Note that the proof 
for the general case is really the same as the one we give in 
Sect.~\ref{Rec}, because $G({\cal G},{\cal A},M)(\gamma)\cong U(1)$ 
for any loop $\gamma$. For completeness, we should remark that we have not 
proved 
that $C_2^2(M,*)(c_{*})\cong\pi_2^2(M,*)$, because of the little technical 
gap that we explained above. Hypothetically $C_2^2(M,*)(c_{*})$ 
might only be a 
``thin quotient'' of $\pi_2^2(M,*)$. However, this technicality is of no 
importance for our approach and we can work with $\pi_2^2(M,*)$ in 
Sect.~\ref{Rec} without any difficulty.

\label{nonsimp}
\section{Barrett's lemma for 2-loops}
This section is a short intermezzo with two technical lemmas. The first one 
we needed in Sect.~\ref{HLB} for the reconstruction of the connection in 
a line-bundle from its holonomy, and the second lemma we will need 
for our reconstruction of the 1-connection in the gerbe obtained from a 
2-holonomy in Sect~\ref{Rec}. In~\cite{Ba91} Barrett proved 
the following lemma, which henceforth we call 
{\it Barrett's lemma} (strictly speaking this is Caetano and Picken's~\cite{CP94} 
version, but the proofs are the same). We state and prove the theorem for 
the case $G=U(1)$, which we need here, but it is true for any Lie group $G$.
\begin{Lem}\label{BLem} The trivial loop extremizes every holonomy 
$\H\colon\pi^1_1(M)\to U(1)$, i.e., given any smooth family of loops 
$\psi\colon[0,1]\to\Omega^{\infty}(M)$ such that $\psi(0)$ is the trivial 
loop, we have $$\frac{d\,\H\circ\psi}{ds}(0)=0,$$ for any holonomy $\H$. Here 
we denote by $\psi$ also the composite of $\psi$ with the natural projection 
$\Omega^{\infty}(M)\to\pi_1^1(M)$.
\end{Lem}
{\bf Proof}: In a neighborhood of $s=0$ all loops in the family are contained 
in one coordinate chart, so it suffices to consider the case $M=\R^n$ with 
basepoint $0\in\R^n$. Using the canonical coordinates in $\R^n$ we can 
write 
$$\psi(s)(t)=(\psi_1(s,t),\psi_2(s,t),\ldots,\psi_n(s,t)).$$ 
This leads to the smooth function $\phi\colon[0,1]^n\times[0,1]\to\R^n$ 
defined by 
$$\phi(s^1,s^2,\ldots,s^n,t)=
(\psi_1(s^1,t),\psi_2(s^2,t),\ldots,\psi_n(s^n,t)).$$ Note that 
$\phi$ defines a smooth family of loops.
We can now write $\psi=\phi\circ\Delta$, where $\Delta\colon[0,1]\to[0,1]^n$ 
is the diagonal map $\Delta(s)=(s,s,\ldots,s)$. A short calculation gives 
$$\frac{d\,\H\circ\psi}{ds}(0)=\frac{d\,\H\circ\phi\circ\Delta}{ds}(0)=
\sum_{i=1}^{n}\frac{\partial\,\H\circ\phi}{\partial s^i}(0,0,\ldots,0)=0.$$
The last equality is a consequence of the fact that all partial derivatives 
are equal to zero. Let us show this for the first partial derivative. 
The value of $\partial\,\H\circ\phi/\partial s^1(0,0,\ldots,0)$ 
only depends on 
the behaviour of $\H\circ\phi$ on the first axis. On this axis we have 
$\phi(s^1,0,\ldots,0)(t)=(\psi_1(s^1,t),0\ldots,0)$. It is not hard to see 
that this is thin homotopic to the trivial loop: for example, a thin homotopy 
is given by $H(s,t)=(\beta(s)\psi_1(s^1,t),0,\ldots,0)$, where $\beta$ is 
the function defined in Sect.~\ref{HLB}. 
Therefore, $\H\circ\phi$ is constant on the first 
axis and its first partial derivative is zero.\qed
\noindent There is a subtlety to be noted here: Lemma~\ref{BLem} is only true when 
$\psi(0)$ is equal to the trivial loop and not when it is thin homotopic to 
it. This is rather important to note, because otherwise one might get the 
wrong impression that the connection constructed in \cite{Ba91,CP94} 
and the gerbe-connection in Sect.~\ref{Rec} of this paper vanish identically 
always.  

In this paper we need the following analogue of Lem.~\ref{BLem} for 2-loops and 2-holonomies (see Def.~\ref{2hol}).
\begin{Lem}\label{2BL} Let $\psi\colon[0,1]^2\to\Omega^{\infty}_2(M)$ be any 
smooth 2-parameter family of 2-loops such that $\psi(0,0)$ is the trivial 
2-loop. Then we have 
$$\frac{\partial^2\,\H\circ\psi}{\partial r\partial s}(0,0)=0,$$ for any 
2-holonomy $\H\colon\pi_2^2(M)\to U(1)$.
\end{Lem}
{\bf Proof}: The proof follows the line of reasoning in Barrett's lemma. Again 
it suffices to consider the case $M=\R^n$ with basepoint $0\in\R^n$. With 
respect to the standard coordinates in $\R^n$ we can write
$$\psi(r,s)(t^1,t^2)=(\psi_1(r,s,t^1,t^2),\psi_2(r,s,t^1,t^2),\ldots,
\psi_n(r,s,t^1,t^2))$$
and
$$\phi(r^1,s^1;r^2,s^2;\ldots,r^n,s^n;t^1,t^2)=$$
$$(\psi_1(r^1,s^1;t^1,t^2),
\psi_2(r^2,s^2;t^1,t^2),\ldots,\psi_n(r^n,s^n;t^1,t^2)).$$
The latter defines a smooth function from $[0,1]^{2n}\times[0,1]^2$ to $\R^n$. 
Note that $\phi$ defines a smooth family of 2-loops. Using the diagonal function $\Delta\colon[0,1]^2\to[0,1]^{2n}$, defined by 
$\Delta(r,s)=(r,s;r,s;\ldots;r,s)$, we can write 
$\psi=\phi\circ\Delta$. Again a short calculation gives 
$$\frac{\partial^2\,\H\circ\psi}{\partial r\partial s}(0,0)=
\frac{\partial^2\,\H\circ\phi\circ\Delta}
{\partial r\partial s}(0,0)=\sum_{i,j=1}^{n}
\frac{\partial^2\,\H\circ\phi}{\partial r^i\partial s^j}(0,0;\ldots;0,0)=0.$$
The last equation is a consequence of the fact that all second order 
partial derivatives 
are equal to zero. Let us show this for the case $i=1, j=2$. In this case 
the value of $\partial^2\,\H\circ\phi/\partial r^1\partial s^2(0,0;
\ldots;0,0)$ 
only depends on the behaviour of $\H\circ\phi$ on the plane spanned by the 
axes corresponding to $r^1$ and $s^2$. In this plane we have 
$$\phi(r^1,0;0,s^2;0,0;\ldots;0,0;t^1,t^2)=(\psi_1(r^1,0;t^1,t^2),
\psi_2(0,s^2;t^1,t^2),0,\ldots,0).$$ 
Now $H(s,t^1,t^2)=(\beta(s)\psi_1(r^1,0;t^1,t^2),\beta(s)
\psi_2(0,s^2;t^1,t^2),0,\ldots,0)$ defines a thin homotopy between the latter 
and the trivial 2-loop, 
so $\H\circ\phi$ is constant on our plane, whence  
$\partial^2\,\H\circ\phi/\partial r^1\partial s^2(0,0)=0$.\qed  
\noindent Again let us note that $\psi(0,0)$ in the previous lemma has to be equal to 
the trivial 2-loop and not just thin homotopic to it.

\label{BL}
\section{The 1-connected case}
Let $M$ be 1-connected.  
\begin{Thm}\label{Recon} Given an arbitrary 2-holonomy 
$\H\colon\pi^2_2(M)\to U(1)$, 
there exists a gerbe $\G$ with gerbe-connection $\A$ on $M$ such that 
the holonomy of $(\G,\A)$ is equal to $\H$. This construction establishes a bijective correspondence 
between equivalence classes of gerbes with gerbe-connections and 2-holonomies.
\end{Thm}
{\bf Proof}: We prove this theorem in several parts. First we show how to 
construct $\G$ and the 0-connection, i.e., the transition 
line-bundles with connections on double intersections and the covariantly 
constant sections on triple intersections. After that we show how to 
construct the 1-connection, $\A^1$, in $\left(\G,\A^0\right)$. In the final 
part we prove the last claim in Thm.~\ref{Recon}.

\noindent{\bf Part 1}: Let $\H\colon\pi^2_2(M)\to U(1)$ be an arbitrary 
2-holonomy. We assume that the covering $\left\{U_i, i\in J\right\}$ of $M$ 
is such that for every $i\in J$ there is a diffeomorphism 
$\phi_i\colon U_i\to B(0,1)\subset\R^n$, where $B(0,1)$ is the unit ball in 
$\R^n$, and that for every pair $i,j\in J$ there is a diffeomorphism 
$\phi_{ij}\colon U_{ij}\to B(0,1)\subset\R^n$ as well. We denote by $x_i$ and 
$x_{ij}$ the centers of $U_i$ and $U_{ij}$ respectively, i.e., 
$x_i=\phi_i^{-1}(0)$ and $x_{ij}=\phi_{ij}^{-1}(0)$. We also assume that 
every point $x\in M$ is the center of some open set in the covering, which we 
denote by $U_x$ when needed. Finally we assume that all $n$-fold 
intersections are contractible. Since $M$ is 1-connected, 
we can choose a path from $*$, the base point in $M$, to $x_i$ for 
every $i\in J$, and a path from $*$ to $x_{ij}$ for every pair 
$i,j\in J$ such that $U_{ij}\ne\emptyset$. For every $i\in J$ we can define a canonical path, 
from $x_i$ to any other point $x\in U_i$ by $\phi_i^{-1}(r_x)$, where 
$r_x$ is the straight line, the segment of a ray, 
in $B(0,1)$ from $\phi_i(x_i)=0$ to $\phi_i(x)$. In particular there is a 
canonical path from $x_i$ to $x_{ij}$, for every $i,j\in J$. Similarly we 
define the canonical path from $x_{ij}$ to any point $y\in U_{ij}$. Having 
chosen all these paths we now have to choose some homotopies. 
Note that the path, denoting the chosen paths from above by arrows, 
$$*\rightarrow x_i\rightarrow x_{ij}\rightarrow *$$ 
is homotopic to $*$, the trivial 
loop, because $M$ is simply-connected, so we can choose a homotopy, 
$P_{ij}$, between them (starting at $*$). By convention $P_{ji}$ is the 
analogous homotopy for the loop 
$$*\rightarrow x_j\rightarrow x_{ij}\rightarrow *.$$
We are now ready to start our construction. 

\begin{figure}
\centerline{
\epsfbox{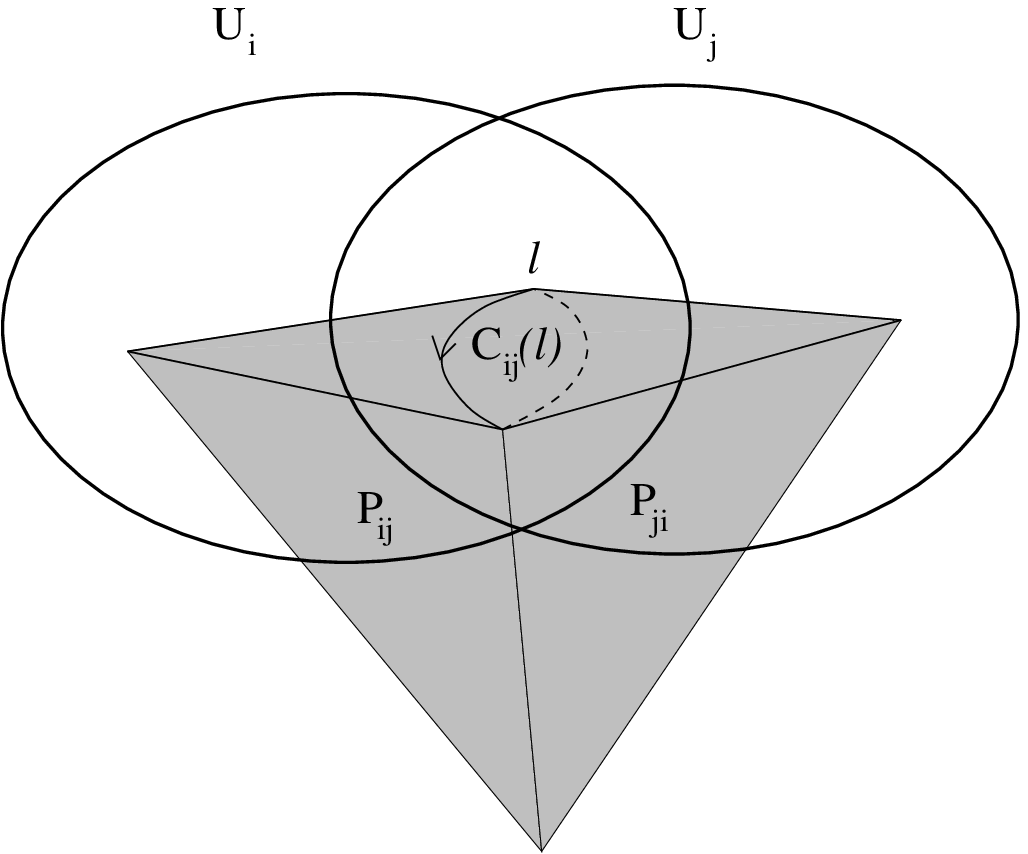}
}
\caption{$s_{ij}(\ell)$}
\label{reccone}
\end{figure}

Choose any pair $i,j\in J$. Let $\ell$ be a loop in $U_{ij}$ based at 
$x_{ij}$. 
Consider the loop $\phi_{i}\circ \ell$ in $B(0,1)$. We can now define 
the cone on $\phi_i\circ \ell$ 
in $B(0,1)$ with top vertex $0\in B(0,1)$. This is just the homotopy between 
$0\in B(0,1)$ and $\phi_i\circ \ell$ obtained by taking all rays from $0$ to 
any point on $\phi_i\circ \ell$ together. Now take the image of this 
cone in $U_i$. We obtain an analogous cone in $U_j$. Next we glue one 
cone onto the other, which corresponds to composing one 
homotopy with the inverse of the other, to obtain a double cone 
$C_{ij}(\ell)$. 
Finally we recenter $C_{ij}(\ell)$ by using $P_{ij}$ and $P_{ji}$ to 
obtain a 2-loop based at $*$, which we denote by $s_{ij}(\ell)$. See 
Fig.~\ref{reccone} for a graphical explanation of our construction. 
Two things are 
immediately clear from this construction: if $\ell$ and $\ell'$ are thin 
homotopic, 
then $s_{ij}(\ell)$ and $s_{ij}(\ell')$ are thin homotopic as well. This holds 
true, because $\mbox{rk}D(CH)$, the rank of the differential of the cone on 
a homotopy $H$, is at most $\mbox{rk}DH+1$. Another obvious observation is 
that $s_{ij}(\ell\star\ell') \stackrel{2}{\sim} s_{ij}(\ell)\star s_{ij}
(\ell')$ 
for any two loops $\ell$ and $\ell'$. Thus we can define a holonomy 
$\H_{ij}\colon\pi_1^1(U_{ij})\to U(1)$ by 
$$\H_{ij}(\ell)=\H(s_{ij}(\ell)).$$
By the results explained in Sect.~\ref{HLB} we obtain a line-bundle 
$\Lambda_{ij}$ with connection $\nabla_{ij}$ and curvature 
$K(\nabla_{ij})$ on $U_{ij}$. By construction $\Lambda_{ji}\cong
(\Lambda_{ij})^{-1}$. Note also that our assumption that intersections are 
contractible implies that $\Lambda_{ij}$ is equivalent to the trivial 
line-bundle. 
Choose a nowhere zero section $\sigma_{ij}$ in $\Lambda_{ij}$.  
On a triple intersection $U_{ijk}$ the tensor product 
$\Lambda_{ijk}=
\Lambda_{ij}\otimes\Lambda_{jk}\otimes\Lambda_{ki}$ is trivial as well and 
$\sigma_{ijk}=\sigma_{ij}\otimes\sigma_{jk}\otimes\sigma_{ki}$ defines a 
nowhere zero section.
Let $\nabla_{ijk}$ denote the connection on $\Lambda_{ijk}$ induced by 
$\nabla_{ij},\nabla_{jk},\nabla_{ki}$ (note that $\nabla_{ji}= -\nabla_{ij}$).
By the results in Sect.~\ref{HLB} we know that the holonomy of $\Lambda_{ij},
\nabla_{ij}$ is exactly equal to $\H_{ij}$, so we conclude that the holonomy 
of $\left(\Lambda_{ijk},\nabla_{ijk}\right)$ around any loop in $U_{ijk}$ is 
trivial, 
because in the construction above we go around each cone twice in opposite 
directions. This means that $\nabla_{ijk}$ is flat. 

\begin{figure}
\vbox{\hskip1cm\epsfbox{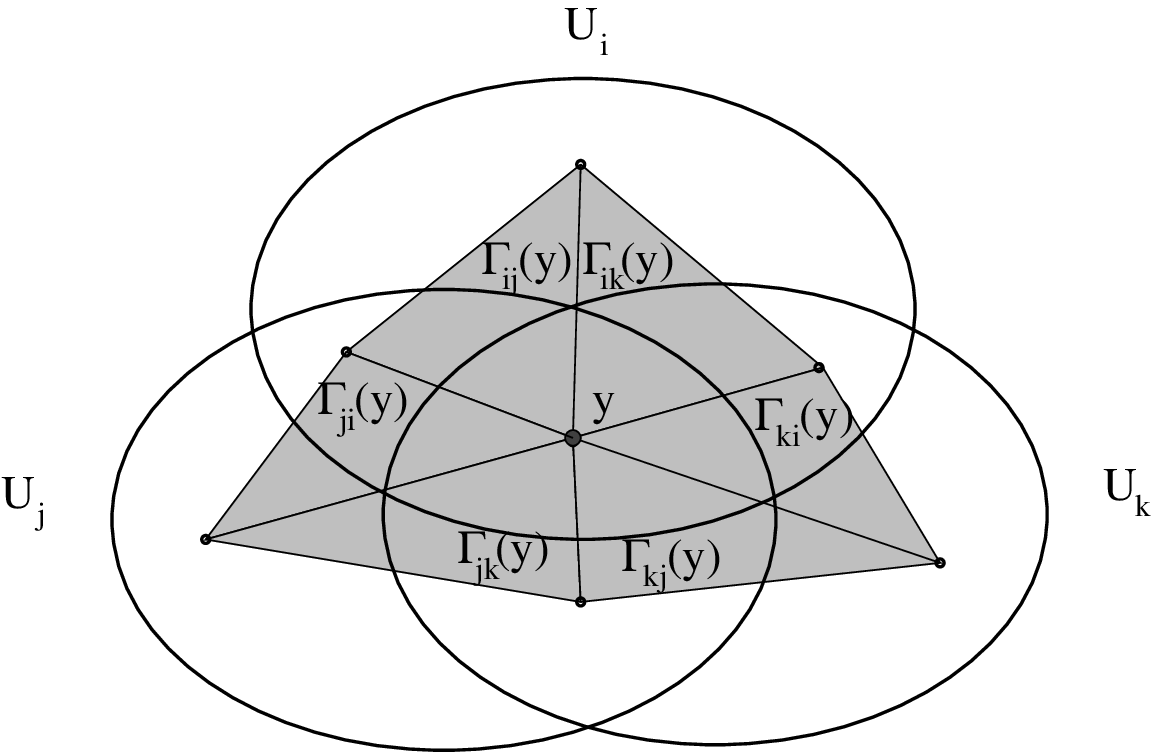}
}
\caption{$s_{ijk}(y)$}
\label{horsect}
\end{figure}

\begin{figure}
\centerline{
\epsfbox{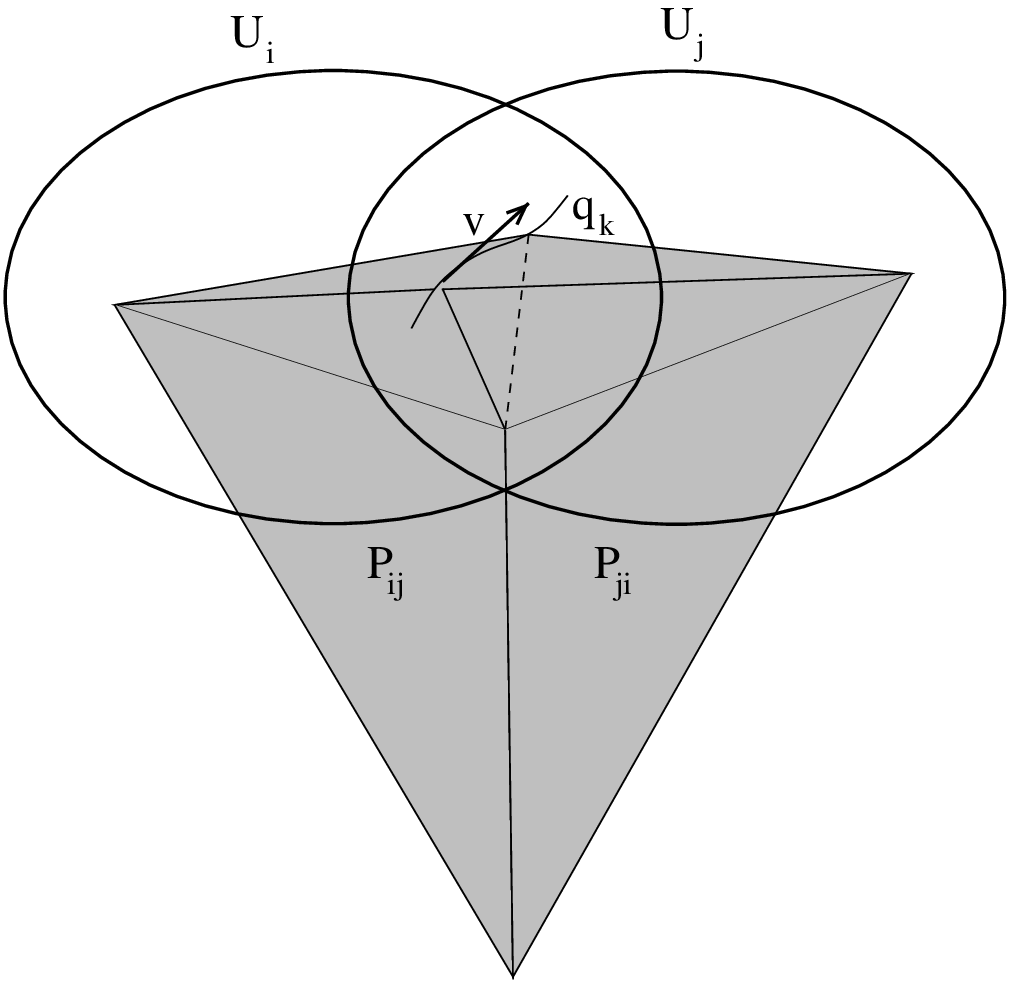}
}
\caption{$s_{ij}(k)$}
\label{0conn}
\end{figure}

\begin{figure}
\centerline{
\epsfbox{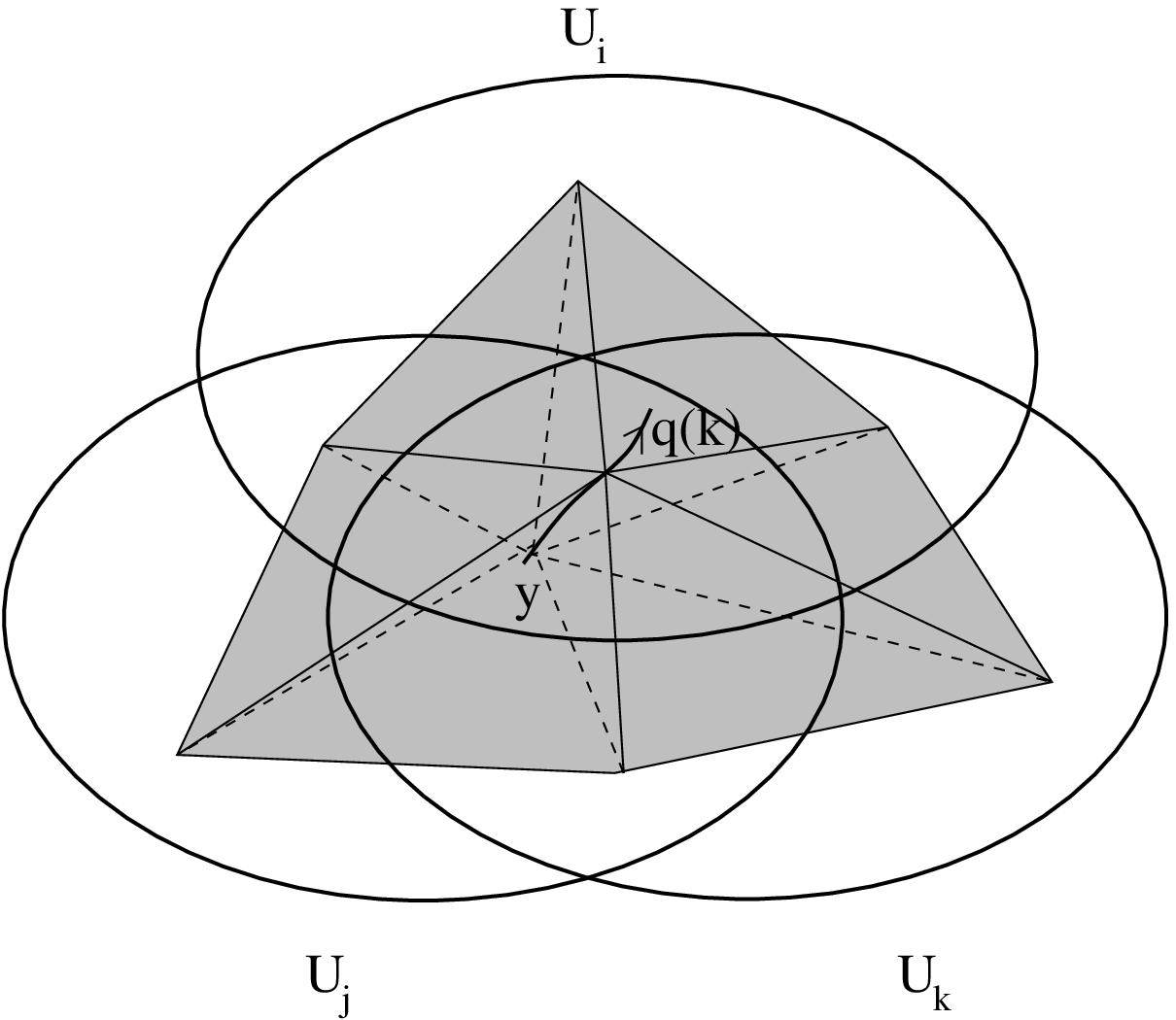}
}
\caption{Proof of Eq.~\ref{eq:0conn}}
\label{0proof}
\end{figure}

Let us now define the 
desired horizontal section $\theta_{ijk}$ in $\Gamma(U_{ijk},\Lambda_{ijk})$. 
For a given point $y\in U_{ijk}$ we define in Fig.~\ref{horsect} a 2-loop 
$s_{ijk}(y)$. 
For example, the part of $s_{ijk}(y)$ relating to the patches $U_i$ and 
$U_j$ is defined 
by the composite of the homotopies $P_{ij}$, $\Gamma_{ij}(y)$, and the 
inverses of $\Gamma_{ji}(y)$ and $P_{ji}$, where $\Gamma_{ij}(y)$ is one of 
the homotopies drawn in Fig.~\ref{horsect}. We now define 
the function $g_{ijk}\colon U_{ijk}\to U(1)$ by 
\begin{equation}
\label{g}
g_{ijk}(y)=\H(s_{ijk}(y)).
\end{equation} 
By construction we have $g_{p(i)p(j)p(k)}=g_{ijk}^{\epsilon(p)}$ for 
any permutation $p\in S_3$, where $\epsilon(p)$ is the sign of $p$. 
It is also easy to see that the collection of functions $g=\left\{g_{ijk}\ 
\vert\ i,j,k\in J\right\}$ defines a \v{C}ech cocycle, i.e., 
$\delta g\equiv 1$. We define $\theta_{ijk}=g_{ijk}\sigma_{ijk}$. The cocycle 
condition satisfied by $g$ implies that $\delta\theta\equiv 1$, because 
$\delta\sigma$ is isomorphic to the canonical section in the trivial line 
bundle by definition. We also claim that for each triple $i,j,k\in J$ the 
section $g_{ijk}\sigma_{ijk}\in\Gamma(U_{ijk},\Lambda_{ijk})$ is covariantly 
constant with respect to $\nabla_{ijk}$. In order to see why this holds true 
we first have to know what 1-form $A_{ij}\in\Omega^1(U_{ij})$ corresponds 
to $\nabla_{ij}$ (remember that we have chosen a section 
$\sigma_{ij}\in\Gamma(U_{ij},\Lambda_{ij})$ which we can use to pull back the 
connection 1-form on the bundle to a 1-form on $U_{ij}$). The results in 
Sect.~\ref{HLB} and our construction of $\left(\Lambda_{ij},\nabla_{ij}
\right)$ show that 
we can define $A_{ij}$ in the following way: let $v$ be any vector in 
$T_y(U_{ij})$, where $y\in U_{ij}$ is an arbitrary point. Represent $v$ by a 
curve $q\colon I\to U_{ij}$, such that $q(0)=y$ and $\dot{q}(0)=v$. Let 
$q_k$ be defined as in Sect.~\ref{HLB}. We can now form the loop $\ell(k)$
$$x_{ij}\rightarrow y\stackrel{q_k}{\rightarrow}q_k(1)=q(k)
\rightarrow x_{ij}.$$
Just as in the beginning of this section we can take the cones on $\ell(k)$ in $U_i$ and $U_j$ respectively and glue them together to form the 
2-loop $s_{ij}(k)$ 
(see Fig.~\ref{0conn}). 
The 
results in Sect.~\ref{HLB} now show that we have 
$$A_{ij}(v)=-i\frac{d}{dk}\log\H(s_{ij}(k)\vert_{k=0}.$$
In order to show that $\nabla_{ijk}(g_{ijk}\sigma_{ijk})=0$ we now have to 
prove the equation 
\begin{equation}
i(A_{ij}-A_{ik}+A_{jk})=-d\log g_{ijk}=
-g_{ijk}^{-1}dg_{ijk}.
\label{eq:0conn}
\end{equation}
Choose a point $y\in U_{ijk}$, a vector $v\in T_y(U_{ijk})$, and a curve 
$q\colon I\to U_{ijk}$ representing $v$. Fig.~\ref{0proof} shows that the 
2-loops defined by $s_{ij}(k)\star s_{ik}(k)^{-1}\star s_{jk}
(k)$ 
and $s_{ijk}(y)\star s_{ijk}(q(k))^{-1}$ are thin homotopic, so the 
holonomy $\H$ maps them to the same number. This proves the desired equation 
after taking derivatives of the logarithms.

\begin{figure}
\centerline{
\epsfbox{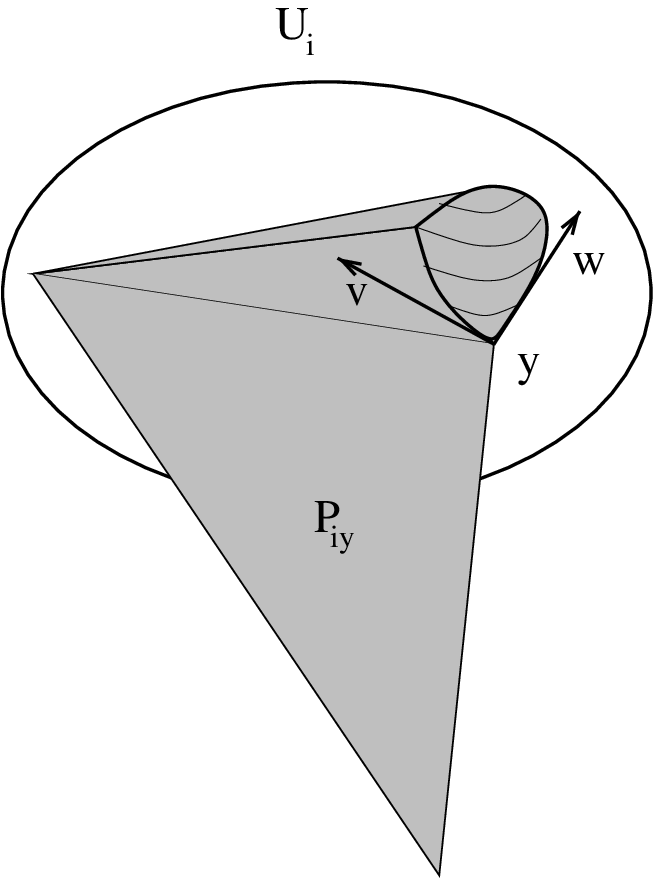}
}
\caption{1-connection}
\label{1connec}
\end{figure}

\begin{figure}
\centerline{
\epsfbox{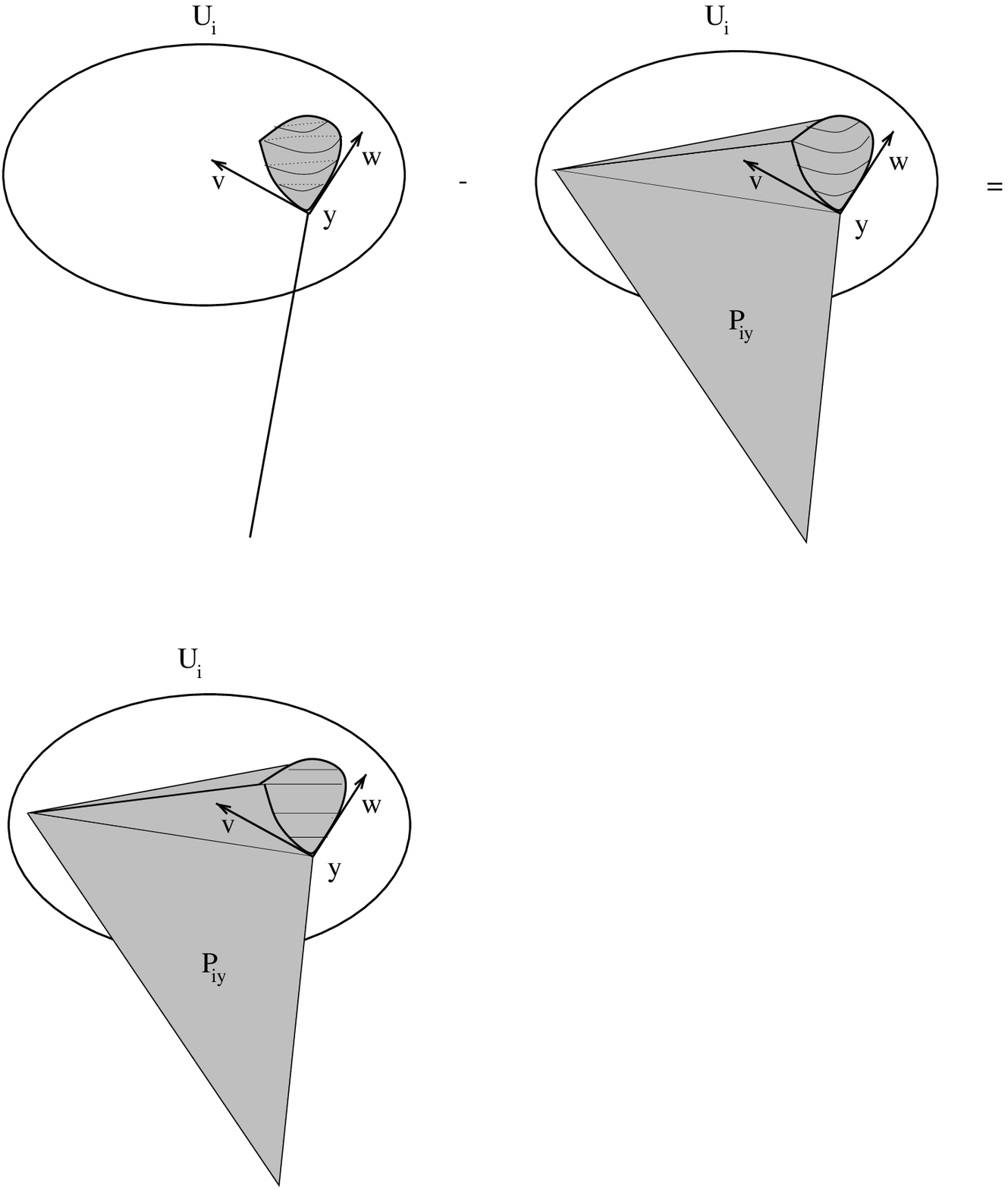}
}
\caption{1-connection is well defined}
\label{1connpro}
\end{figure}

\noindent{\bf Part 2}: In this part we are going to define the 2-forms 
$F_i\in U_i$ 
which constitute the 1-connection on $\G$. Let $y\in U_i$ be an arbitrary 
point, and $v,w\in T_y(U_i)$ two (linearly independent) vectors. In a small 
neighborhood $W\subset U_i$ of $y$, we 
can choose two commuting flows $q,r\colon I\times W\to U_i$, representing 
$v$ and $w$ respectively, i.e., $q(0,y)=r(0,y)=y$, $\dot{q}(0,y)=v$, 
$\dot{r}(0,y)=w$, where $\dot{q}$ means the time derivative 
(first coordinate) of $q$. We say that $q$ and $r$ commute when 
$q(t_1,r(t_2,x))=r(t_2,q(t_1,x))$ for all $t_1,t_2\in I$ and all $x\in W$ 
(this corresponds exactly to saying that the two vector fields induced by 
$q$ and $r$ commute). Locally we can always choose such flows, because it 
is possible in $\R^n$. In particular, for any $k,l\in I$, we can define the 
2-path, starting in $y$, by 
$$r\left(\beta\left(t_2\right)l,q\left(\beta\left(t_1\right)k,y\right)
\right).$$ 
See Sect.~\ref{HLB} for the 
definition of $\beta$. Now take the cone on the boundary of this 2-path 
with vertex in $x_i$ and glue this cone on top of the 2-path in order to get 
a 2-loop, centered at $y$. As before we have to recenter this 
2-loop. One of our assumptions was that every point $x$ is the center of some 
open $U_x$, so we recenter using $P_{iy}$ to obtain a 2-loop 
based at $*$ (see Fig.~\ref{1connec}), which we denote by $s_i(k,l)$. Now 
define $F_i\in\Omega^2(U_i)$ by 
\begin{equation}
F_i(v,w)=-i\frac{\partial^2}{\partial k\partial l}\log\H(s_i(k,l))
\vert_{(k,l)=(0,0)}.
\label{eq:1conn}
\end{equation}   
We have to show that $F_i(v,w)$ is well defined, i.e. independent of the 
choice of flows, and that the set $\left\{F_i\ \vert\ i\in J\right\}$ defines 
a 1-connection. Both facts are consequences of the same observations, 
which we explain now. Let $y,v,w,q,r$ be as above. First of all we claim that 
$F_y(v,w)=0$ for any $q,r$ representing $v,w$. Here $F_y$ is defined using 
the open set $U_y$ whose center is $y$. This follows from our next order version 
of Barrett's lemma, which is Lem.~\ref{2BL}. Note that, in the notation from 
above, if we recenter our 2-loops so that they become based at $y$, 
the smooth 2-parameter family of 2-loops $s_y(k,l)$, depending on the 
parameters $k,l\in I$, satisfies the condition of Lem.~\ref{2BL} 
because $s_y(0,0)=y$, the constant 2-loop at $y$. Clearly this is 
true for 
any flows representing $v$ and $w$. Next, let us have a look at 
Fig.~\ref{1connpro}. 
In this figure we show that the 2-loops 
$s_y(k,l)\star s_i(k,l)^{-1}$ and $s_{iy}(k,l)$ are 
thin homotopic. Here $s_{iy}(k,l)$ is defined analogously to $s_{iy}(k)$ 
(see Fig.~\ref{0conn}) using the loop $\ell(k,l)$ around the 
boundary of the 2-path 
which we used in the 
definition of $F_i$.  
Therefore we have 
\begin{equation}
\H(s_y(k,l))\H(s_i(k,l))^{-1}=\H(s_{iy}(k,l)).
\label{eq:HHH}
\end{equation}
Now, the right-hand side of this equation is exactly the holonomy of 
$\nabla_{iy}$ around the loop $\ell(k,l)$. Applying Stokes' theorem to the 
pull-back, $F_{iy}=dA_{iy}\in\Omega^2(U_{iy})$, of the curvature 2-form 
$K(\nabla_{iy})$ via the section $\sigma_{iy}$, and taking the second order 
partial 
derivative on the right-hand side 
gives  
$$-i\frac{\partial^2}{\partial k\partial l}
\log\H(s_{iy}(k,l))\vert_{(k,l)=(0,0)}=F_{iy}(v,w).$$
Taking also the corresponding second order partial derivative on the 
left-hand side of eq.~\ref{eq:HHH} gives us 
\begin{equation}
F_y(v,w)-F_i(v,w)=F_{iy}(v,w).\label{eq:FFF}
\end{equation}
This equation shows two things at once. In the first place we conclude that 
$F_i$ is well defined, because we have 
$$F_i(v,w)=F_i(v,w)-F_y(v,w),$$ 
since the last term is zero, and $F_{iy}(v,w)$ does not depend on the choice of flows, 
because $F_{iy}$ is an honest 2-form. Secondly eq.~\ref{eq:FFF} implies that 
the 
$F_i$ define a 1-connection in $\G$, because for any $j\in J$ we now get 
$$
\begin{array}{rcl}
F_j(v,w)-F_i(v,w) & = & F_j(v,w)-F_y(v,w)+F_y(v,w)-F_i(v,w) \\
                  & = & F_{yj}(v,w)+F_{iy}(v,w) \\
                  & = & F_{ij}(v,w).
\end{array}
$$ 
The last equality follows from the fact that the $F_{ij}$ are curvature 
2-forms of the connections $A_{ij}$ which define a 0-connection on $\G$, i.e., 
$\delta A=-id\log g$.    
 
\begin{figure}
\centerline{
\epsfbox{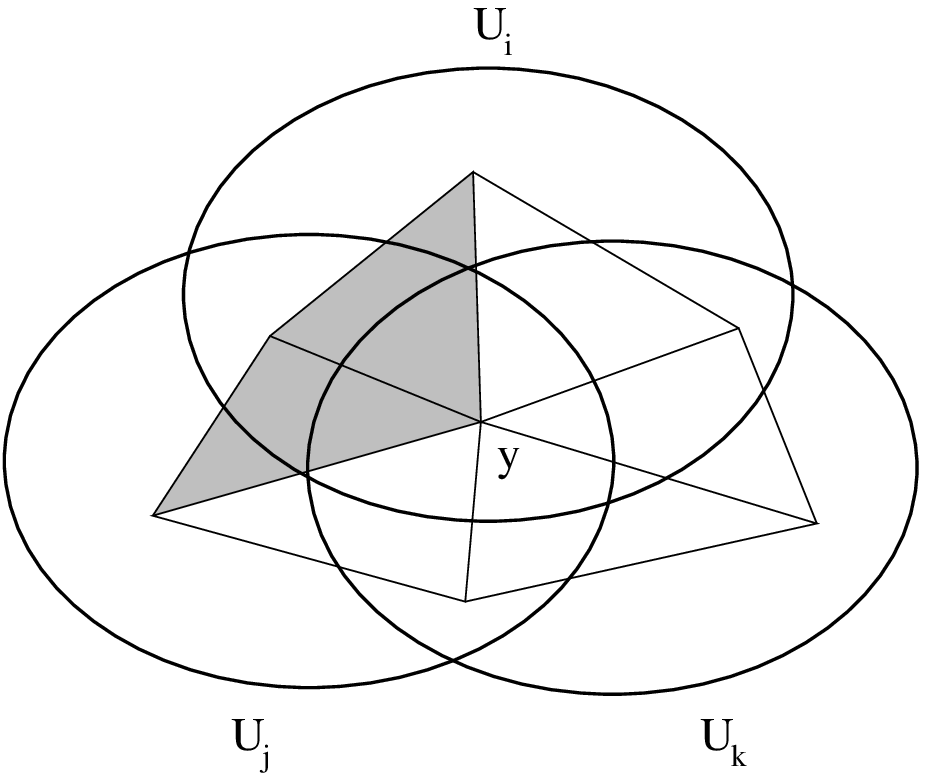}
}
\caption{$h_{ij}$}
\label{0ga}
\end{figure}

\begin{figure}
\centerline{
\epsfbox{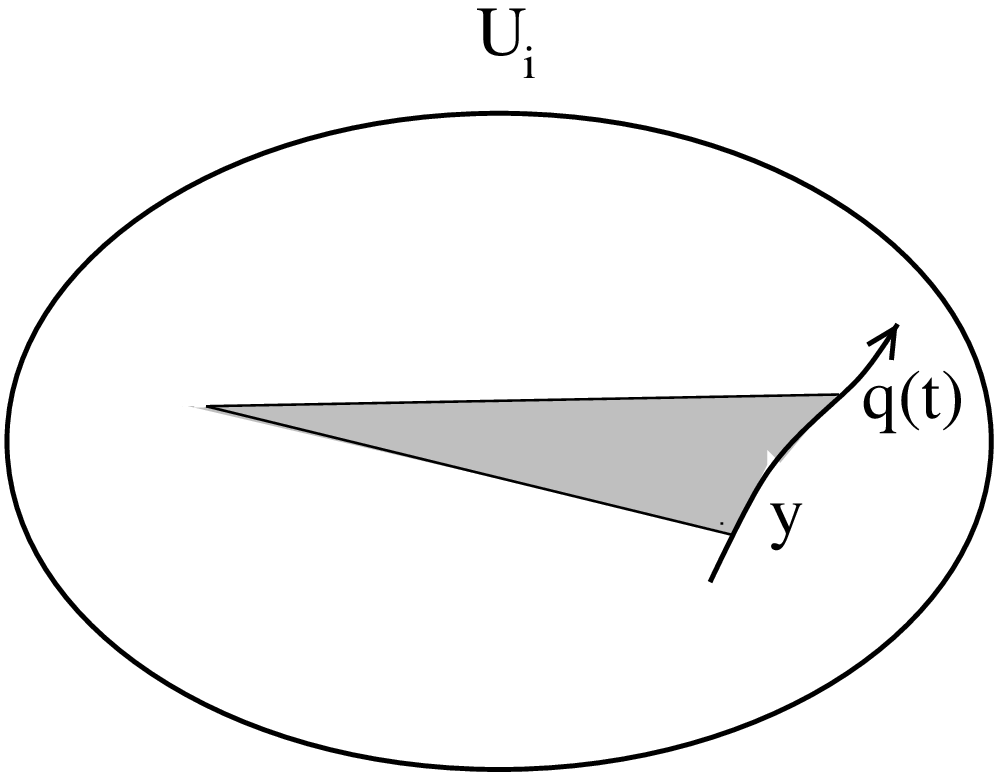}
}
\caption{$B_i$}
\label{1ga}
\end{figure}

\noindent{\bf Part 3}: In this final part of the proof of Thm.~\ref{Recon} 
we show that the construction above defines a bijection between equivalence 
classes of gerbes with gerbe-connections on the one hand and 
2-holonomies on the other. Let $\G$ be a 
gerbe on $M$ and $\A$ a gerbe-connection in $\G$, and let 
$\H\colon\pi_2^2(M)\to U(1)$ be the gerbe-holonomy of $\G,\A$. Using 
the construction above we obtain a new gerbe $\G'$ with a gerbe-connection 
$\A'$ from $\H$. Let us show that $\G,\A$ and $\G',\A'$ are 
equivalent. Our proof is local, so we assume that $\G$ (resp. $\G'$) is given 
by a cocycle $g_{ijk}$ (resp. $g'_{ijk}$) and that $\A$ (resp. $\A'$) is 
given by $A_{ij},F_i$ (resp. $A'_{ij},F'_i$). Let $y\in U_{ijk}$ be an 
arbitrary point. Recall (eq.~\ref{g}) how we defined $g'_{ijk}(y)\in U(1)$: 
$$g'_{ijk}(y)=\H(s_{ijk}(y)).$$
At the end of Sect.~\ref{HG} we obtained a concrete formula for 
$\H(s)$, for any $s\in\pi_2^2(M)$. We define the function 
$h_{ij}\colon U_{ij}\to U(1)$ by 
\begin{equation}
h_{ij}(y)=\exp\left(i\int_{I^2_{ij}}\epsilon_s\right)\cdot
\lambda^{-1}_{ij}(y),\label{eq:h}
\end{equation}
where $I^2_{ij}$ is the part of $I^2$ which is mapped onto the part of 
$s_{ijk}(y)$ which goes from $*$ to $U_i$ and $U_j$ (see Fig.~\ref{0ga}) and 
$\lambda_{ij}$ is defined in Sect.~\ref{HG}. 
Likewise we define $h_{ik}$ and 
$h_{jk}$. From formula~(\ref{eq:h}) it is immediately clear that we have 
$h_{ji}=h^{-1}_{ij}$ and   
$$g'_{ijk}(y)=g_{ijk}(y)h_{ij}(y)h_{jk}(y)h_{ki}(y).$$
This shows that $\G$ and $\G'$ are equivalent as gerbes. In order to establish the full equivalence between $\A$ and $\A'$ we also define the 1-forms 
$B_i\in\Omega^{1}(U_i)$ by
$$B_i(v)=-\frac{d}{dk}\left(\int_{C_i(k)}F_i\right)\vert_{k=0},$$ 
where $q(t)$ is a curve in $U_i$ representing $v\in T_y(U_i)$ and 
$C_i(k)$ is 
the 2-path in $U_i$ defined by the cone on $q_k$ with vertex $x_i$ 
(see Fig.\ref{1ga}). 
A small and simple calculation in $\R^n$, which we omit, shows that the 
definition of $B_i(v)$ does not depend on the 
choice of $q(t)$. From our construction of $A'_{ij}$ from $\H$, as 
explained in part 1 of this proof (see Fig.~\ref{0conn}), we get the 
following equality:
$$\H(s)=\exp\left(i\left(
\int_{q_k}A_{ij}+\int_{C_i(k)}F_i-\int_{C_j(k)}F_j\right)\right)h_{ij}(q(0))h_{ij}^{-1}(q(k)).$$ 
Here we take $s=s_{ij}(k)$, which we defined in Fig.~\ref{0conn}. 
Taking the 
derivatives at $k=0$ on both sides gives
$$i(A')_{ij}(v)=iA_{ij}(v)+iB_j(v)-iB_i(v)-\left(d\log h_{ij}(v)\right)(v).$$
Finally we have to prove that $F'_i=F_i+dB_i$ for all $i$. Let 
$v,w\in T_y(U_i)$ be two arbitary vectors. Then 
\begin{eqnarray*}
F'_i(v,w)&=&F'_i(v,w)-F'_y(v,w)\\
&=&(F')_{yi}(v,w)\\
&=&d(A')_{yi}(v,w) \\
&=&d\left(A_{yi}+B_i-B_y+id\log h_{yi}\right)(v,w)\\
&=&\left(F_{yi}(v,w)+dB_i-dB_y\right)(v,w)\\
&=&\left(F_i+dB_i\right)(v,w)-
\left(F_y+dB_y\right)(v,w)\\
&=&\left(F_i+dB_i\right)(v,w).
\end{eqnarray*}
The last equality follows from Barrett's lemma for 2-loops (Lem.~\ref{2BL}), 
because $\left(F_y+dB_y\right)(v,w)$, according to our formula for 
$\H$ at the end of Sect.~\ref{HG}, is equal to 
$$-i\frac{\partial^2}{\partial k\partial l}\log\H
(s_y(k,l))\vert_{(k,l)=(0,0)}=0,$$
following the notation in~(\ref{eq:1conn}). The smooth 2-parameter family of 
2-loops 
$s_y(k,l)$ starts at the trivial 2-loop at $y$ after recentering, so we 
can indeed apply Lem.~\ref{2BL}. 

\begin{figure}
\centerline{
\epsfbox{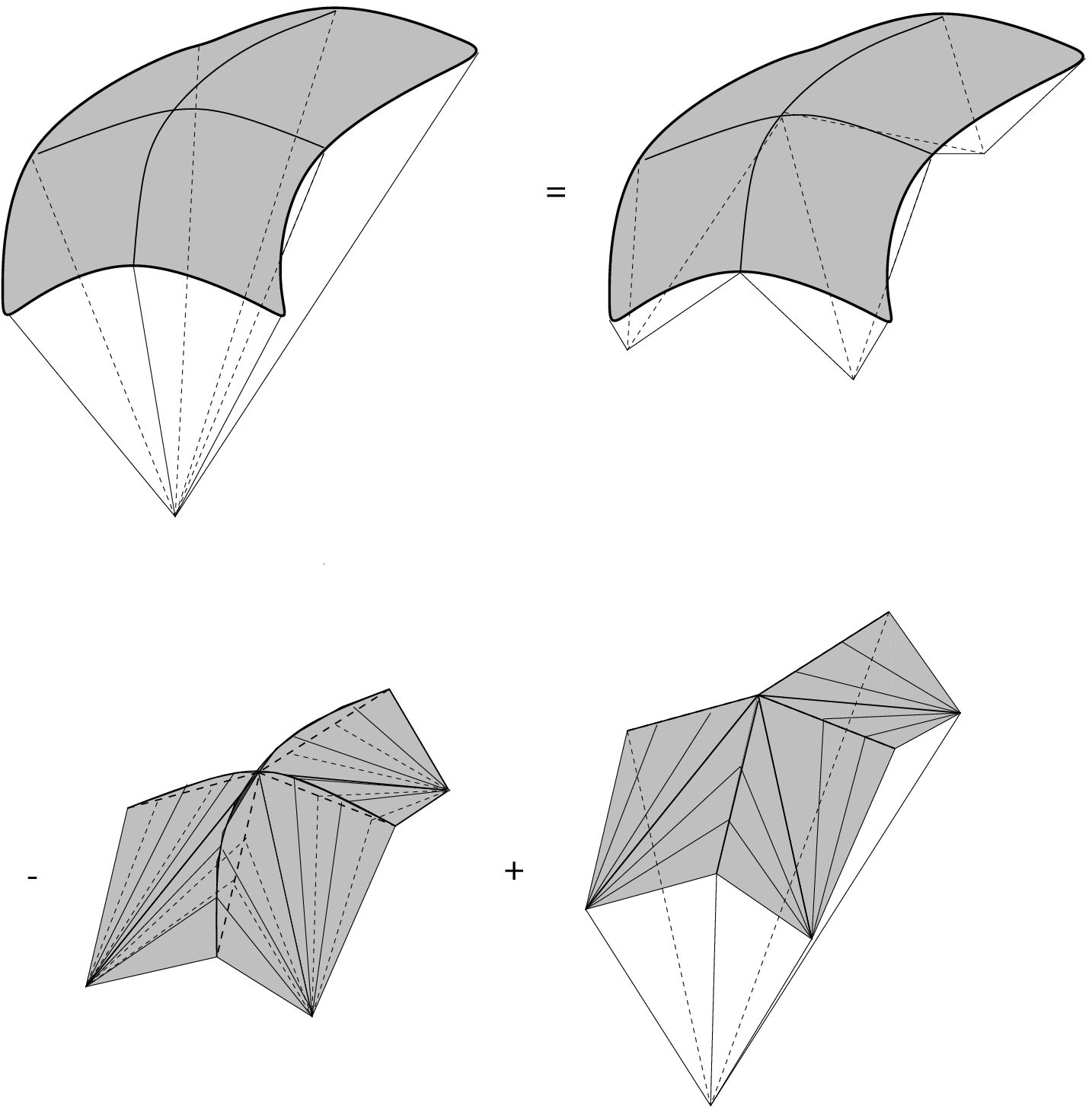}
}
\caption{$\H_{\G,\A}(s)$}
\label{hpro}
\end{figure}

Conversely, let $\H$ be a 2-holonomy, reconstruct 
$\G,\A$ as above, and let $\H_{\G,\A}$ be the 
gerbe-holonomy of $(\G,\A)$. We show that $\H_{\G,\A}=\H$. The analogous 
proof for line-bundles, which we gave in Sect.~\ref{HLB}, relied on the fact 
that 
the holonomy around a loop, $\ell$, can be written as the holonomies around 
many loops each of which only shares a part with $\ell$, such that in the end 
all parts of the loops that do not belong to $\ell$ cancel out. The same idea 
underlies our proof for gerbe-holonomies. Let $s\colon I^2\to M$ be a 
2-loop. In Fig.~\ref{hpro} we have drawn a part of the image of $s$ which is 
covered by four open sets, $U_i,U_j,U_k,U_l$, and which we denote by 
$s_{ijkl}$. 
The formula for $\H_{\G,\A}(s)$ 
at the end of Sect.~\ref{HG} shows that the part of $\H_{\G,\A}(s)$ 
which corresponds to $U_i,U_j,U_k,U_l$ is given by integration of the 2-forms 
$F_i,F_j,F_k,F_l$ over that part of the image of $s$ which intersects $U_i\cup
U_j\cup U_k\cup U_l$, by integration of the 1-forms $A_{ij},A_{jk},A_{kl},
A_{li}$ 
over the edges $E_{ij},E_{jk},E_{kl},E_{li}$ in the image of $s$, and finally 
by evaluating $g_{ijk}g_{ilk}^{-1}$ at $y$, a point in the image of $s$ and 
in the 
intersection $U_i\cap U_j\cap U_k\cap U_l$. The key observation, just as in 
the case 
for line-bundles, is that $F_i,A_{ij},g_{ijk}$ are all defined in terms of 
$\H$. By close inspection of the definition in part 1 and 2 of this proof 
we find that the images of the 2-loops which define 
$F_i,F_j,F_k,F_l$ contain all of $s_{ijkl}$ but they contain more. 
This extra bit 
of 2-loop gets cancelled by the 2-loops defining $A_{ij},A_{jk},A_{kl},
A_{li}$ and 
$g_{ijk}g_{ilk}^{-1}$. In Fig.~\ref{hpro} the 2-loop in the first picture 
corresponds to $s_{ijkl}$. In the second picture we see the 2-loops 
corresponding to $F_i,F_j,F_k,F_l$. After composing with the inverse of 
the 2-loop in the third picture, which corresponds to 
$A_{ij},A_{jk},A_{kl},A_{li}$, we are left with the 2-loop in the 
fourth picture 
which corresponds to $g_{ijk}g_{ilk}^{-1}(y)$. We see 
that all parts of the 2-loops above which are not part of $s_{ijkl}$ cancel 
out and so we are left with $s_{ijkl}$ only. This shows that we have 
$\H_{\G,\H}(s)=\H(s)$.\qed             

\begin{Rem}
\label{transg} 
{\rm In Sect.~\ref{nonsimp} we recalled that a gerbe with 
gerbe-connection 
on $M$ induces a line-bundle with connection on $\Omega(M,x)$, which actually 
quotients to a line-bundle on $\pi_1^1(M,x)$. We remark 
that in the simply-connected case there exists a different construction of 
that line-bundle which 
follows from the results in this section.
 
In~\cite{Ba91,CP94} the authors reconstructed a principal bundle with 
connection from a holonomy by a global method, i.e., they reconstructed the 
total space of the bundle first and showed that there is a natural lifting 
of paths built in which defines a connection. In the case of gerbes one could 
do the same for a given 2-holonomy, as we sketch in the following. 
Let $\H\colon\pi_2^2(M)\to U(1)$ be a 2-holonomy 
($M$ continues to be 1-connected). One can define the total space 
$$P^{\infty}_2(M,*)\times U(1)/\sim.$$
Here $P^{\infty}_2(M,*)$ is the set of all 2-paths $s\colon I^2\to M$ 
such that $s(r,0)=s(r,1)=*\ \forall r\in I$ and $s(0,t)=*,\ \forall t\in I$. 
The equivalence relation is defined by 
$$(s_1,l_1)\sim(s_2,l_2)\Leftrightarrow \forall t\in I\ s_1(1,t)=s_2(1,t)\ 
\wedge\ l_2=\H(s_2\star s_1^{-1})l_1.$$
It is easy to check that this relation indeed defines an equivalence relation. 
The set of equivalence classes is a line-bundle on $\Omega^{\infty}(M,*)$, 
where the projection $\pi\colon P^{\infty}_2(M,*)\times U(1)/\sim\ \to 
\Omega^{\infty}(M,*)$ is defined by 
$$\pi([(s,l)])(t)=s(1,t).$$
One can of course quotient this line-bundle to obtain one over $\pi_1^1(M,x)$, 
following the observations in Sect.~\ref{nonsimp}.
It looks likely that the whole construction carried out 
in~\cite{Ba91,CP94} works in 
this setting as well. For example, the connection would now come in the 
form of a lifting function of paths of loops. However, everything becomes 
infinite-dimensional in such an approach. To avoid that we have opted to do 
everything locally in $M$, which is very concrete although less elegant 
maybe. 
}  
\end{Rem}
\begin{Rem}
{\rm 
In this remark we want to point out a relation between thin homotopy 
groups and hypercohomology groups that is a consequence of our results. 
In ordinary homotopy theory it is 
well known that for an $(n-1)$-connected manifold $M$, the Hurewicz map 
$\pi_n(M)_{ab}\to H_n(M)$ defines an isomorphism of groups. The results in 
\cite{Ba91,CP94} show that for a connected manifold $M$ there exists an isomorphism 
of groups between the group of holonomies $\left\{\pi_1^1(M)\to U(1)\right\}$ 
and the hypercohomology group $H^1(M,\underline{\Co}^*_M\rightarrow
\underline{A}^1_{M,\Co})$. This isomorphism exists because both groups 
classify line-bundles 
with connections up to equivalence. In our case we see 
that for a 1-connected 
manifold $M$, there exists an isomorphism of groups between the group 
of 2-holonomies $\left\{\pi_2^2(M)\to U(1)\right\}$ and the hypercohomology 
group 
$H^2(M,\underline{\Co}^*_M\rightarrow\underline{A}^1_{M,\Co}\rightarrow
\underline{A}^2_{M,\Co})$ (both groups classify 
gerbes with gerbe-connections up to equivalence). It is likely that, in 
general, for an $(n-1)$-connected manifold $M$, 
the groups $\left\{\pi_n^n(M)\to U(1)
\right\}$ and $H^n(M,\underline{\Co}^*_M\rightarrow\underline{A}^1_{M,\Co}
\rightarrow\underline{A}^2_{M,\Co}\rightarrow \cdots\rightarrow
\underline{A}^n_{M,\Co})$ are isomorphic.}
\end{Rem}

\label{Rec}
\vskip0.5cm
\centerline{\bf Acknowledgements}
\vskip.5cm
\noindent The first author thanks John Barrett for helpful discussions and 
Ana Cannas da Silva for indicating Hitchin's paper.

This work was supported by {\em Programa Operacional
``Ci\^{e}ncia, Tecnologia, Inova\c{c}\~{a}o''} (POCTI) of the
{\em Funda\c{c}\~{a}o para a Ci\^{e}ncia e a Tecnologia} (FCT),
cofinanced by the European Community fund FEDER.  
The first author is presently on leave from the Universidade do Algarve 
and working with a postdoctoral fellowship from FCT at the University of 
Nottingham (UK).
\bibliographystyle{plain}

\end{document}